\documentclass[amssymb,12pt]{article}

\usepackage{mathtools}
\usepackage{geometry,amssymb,amsthm,amsmath,
graphics,enumerate}
\usepackage{times}
\usepackage{amsfonts}
\usepackage{amscd}
\usepackage{url}
\usepackage{enumitem}
\author{ Michael C.\
Laskowski  
and Danielle S.\ Ulrich
 \thanks{Both authors partially supported
by NSF grant DMS-2154101.}
\\
Department of Mathematics\\University of Maryland
}
\newbox\smilebox
\newbox\anchorbox
\newbox\noanchorbox
\newbox\tempbox

\setbox\smilebox=\hbox{$\smile$}

\def\anchor{\hbox{\vtop{
           \hbox to \wd\smilebox{\hfil\vrule width.4pt height7pt depth1pt\hfil}
           \vskip  -11.5truept
           \hbox to \wd\smilebox{\hfil$\smile$\hfil}}}}
\setbox\anchorbox=\anchor
\def\noanchor{\hbox{\vtop{
           \hbox to \wd\anchorbox{\hfil\anchor\hfil}
           \vskip -14truept
           \hbox to \wd\anchorbox{\hfil/\hfil}}}}
\setbox\noanchorbox=\noanchor

\def\fg#1#2#3{\setbox\tempbox=\hbox{$\scriptstyle{#2}$}
\ifnum\wd\anchorbox>\wd\tempbox\dimen255=\wd\anchorbox
\else\dimen255=\wd\tempbox\fi
{#1\,\vtop{\hbox to \dimen255{\hfil\anchor\hfil}
           \vskip -6truept
           \hbox to \dimen255{\hfil$\scriptstyle{#2}$\hfil}}
           \,#3}}

\def\nfg#1#2#3{\setbox\tempbox=\hbox{$\scriptstyle{#2}$}
\ifnum\wd\noanchorbox>\wd\tempbox\dimen255=\wd\noanchorbox
\else\dimen255=\wd\tempbox\fi
{#1\,\vtop{\hbox to \dimen255{\hfil\noanchor\hfil}
           \vskip -6truept
           \hbox to \dimen255{\hfil$\scriptstyle{#2}$\hfil}}
           \,#3}}

\setbox1=\hbox{$\bot$}

\def\north#1#2{#1\,
\hbox{$\bot$\llap {\hbox to\wd1 {\hfil $/$\hfil}}}
\,#2}

\def\nao#1#2#3{#1\  \hbox{\vtop{ 
\baselineskip=4pt
\hbox{$\bot$\llap {\hbox to\wd1 {\hfil $/$\hfil}}
\hskip .05em \llap{\hbox{$^{\scriptscriptstyle{a}}$}}}\hbox{$\scriptstyle
{#2}$}}}\, #3}

\def\includeE#1{{\lhook\kern-3.5pt\joinrel\smash{
    \mathop{\longrightarrow}\limits^{#1}}}}

\def\efor/{Example~\ref{E4}}

\def\BL/{Baldwin--Lachlan}
\def\Bu/{Buechler}
\def\Hr/{Hrushovski}
\def\lm/{locally modular}
\def\wm/{weakly minimal}
\def\nm/{non--modular}
\def\ss/{superstable}
\def\ud/{unidimensional}
\def\sm/{strongly minimal}

\def\abar{\overline{a}}

\def\bbar{\overline{b}}
\def\cbar{\overline{c}}
\def\dbar{\overline{d}}
\def\ebar{\overline{e}}
\def\fbar{\overline{f}}

\def\hbar{\overline{h}}

\def\mbar{\overline{m}}
\def\nbar{\overline{n}}

\def\ubar{\overline{u}}

\def\xbar{\overline{x}}
\def\ybar{\overline{y}}
\def\zbar{\overline{z}}

\def\tp{{\rm tp}}

\def\tr/{trivial}
\def\nt/{non--trivial}
\def\st/{strong type}

\def\TV/{Tarski--Vaught}

\def\sc/{sound construction}
\def\ac/{atomic construction}

\def\fal/{functional}
\def\upl/{unique parallel lines}
\def\chp/{categorical in a higher power}

\def\includeE#1{{\lhook\kern-3.5pt\joinrel\smash{
    \mathop{\longrightarrow}\limits^{#1}}}}

\def\efor/{Example~\ref{E4}}

\def\BL/{Baldwin--Lachlan}
\def\Bu/{Buechler}
\def\Hr/{Hrushovski}
\def\lm/{locally modular}
\def\wm/{weakly minimal}
\def\nm/{non--modular}
\def\ss/{superstable}
\def\ud/{unidimensional}
\def\sm/{strongly minimal}

\def\abar{\overline{a}}

\def\bbar{\overline{b}}
\def\cbar{\overline{c}}
\def\dbar{\overline{d}}
\def\ebar{\overline{e}}
\def\fbar{\overline{f}}

\def\hbar{\overline{h}}

\def\mbar{\overline{m}}
\def\nbar{\overline{n}}

\def\ubar{\overline{u}}

\def\xbar{\overline{x}}
\def\ybar{\overline{y}}
\def\zbar{\overline{z}}

\def\tp{{\rm tp}}

\def\tr/{trivial}
\def\nt/{non--trivial}
\def\st/{strong type}

\def\abar{\bar{a}}
\def\bbar{\bar{b}}
\def\cbar{\bar{c}}
\def\dbar{\bar{d}}
\def\ebar{\bar{e}}
\def\ybar{\bar{y}}
\def\phi{\varphi}

\def\C{{\mathfrak  C}}

\def\F{{\cal F}}
\def\FF{{\bf F}}

\def\P{{\cal P}}

\def\Z{{\mathbb Z}}
\def\tp{{\rm tp}}

\def\Fa0{{\FF^a_{\aleph_0}}}

\def\<{\langle}
\def\>{\rangle}

\def\U{{\cal U}}
\newtheorem{Theorem}{Theorem}[section]
\newtheorem{Proposition}[Theorem]{Proposition}
\newtheorem{Definition}[Theorem]{Definition}
\newtheorem{Notation}[Theorem]{Notation}
\newtheorem{Remark}[Theorem]{Remark}

\newtheorem{Lemma}[Theorem]{Lemma}
\newtheorem{Corollary}[Theorem]{Corollary}

\newtheorem{Fact}[Theorem]{Fact}

\def\ss{{\bf s}}

\DeclareMathOperator{\Mod}{Mod}

\def\CC{{\mathcal C}}

\def\KK{{\mathbf K_0}}
\def\K1{{\mathbf K_1}}

\def\F{{\cal F}}
\def\Vbar{{\overline{V}}}
\def\Abar{{\overline{A}}}

\newcommand\myrestriction{\mathord\restriction}
\def\mr#1{\myrestriction_{#1}}
\def\0bar{\overline{0}}
\def\F{{\cal F}}

\def\gammabar{{\overline{\gamma}}}
\def\rhobar{{\overline{\rho}}}
\def\xtilde{{\hat{x}}}
\def\RsubR{{}_R R}
\def\L{{\cal L}}
\def\fhat{{\hat{f}}}

\begin{document}

	\title{Borel completeness of $R$-modules when $R$ fails the DCC on pp-definable subgroups}
	
	\date{\today}

	\maketitle
	
	\begin{abstract}  We prove that for any countable ring $R$ (not necessarily commutative), if the associated left $R$-module $\RsubR$ has a strictly descending sequence of pp-definable subgroups, then the theory $Th(R^{(\omega)})$ of the infinite dimensional direct sum is Borel complete.  From this, we conclude that if $R$ is countable and not left perfect, then the theory of $R$-modules is Borel complete, and we give a full characterization of which countable simple rings have Borel complete theories.   
	One special case is that the complete theory $Th(\Z^{(\omega)})$ is Borel complete, which strengthens the existing proofs of the Borel completeness of TFAB, the theory of torsion free abelian groups.
	
	The proof also introduces, relative to the chosen pp-chain, a proper two-sided ideal $L^R$,
	 and a notion of f.g. hulls which, for countable rings and countable parameter sets in theories satisfying $T=T^{\aleph_0}$ exist and are unique up to isomorphism. These constructions may be of independent interest in the model theory of modules.
\end{abstract}

\section{Introduction}  A long-standing project about  modules is -- given a countable ring $R$ --  determine the complexity of $\CC_R$, the class of all countable left $R$-modules, up to isomorphism.  One invariant is given by the number of isomorphism types occurring in $\CC_R$, but this is very coarse, as there are very few possibilities for this number.  

A more robust measure of complexity,  {\em Borel reducibility,}
was introduced by Friedman and Stanley~\cite{FS}.   It was introduced 
to measure the complexity of 
quotients of Polish spaces by natural equivalence relations.  
In model theory, the natural context is the Polish space $Str_\L$ of all structures with universe
$\omega$ in a countable language $\L$.  
With the natural topology on $Str_\L$,  for any $\L$-theory $T$, its class $Mod(T)$ of models with universe $\omega$
 forms a Borel subset of $Str_\L$ and the equivalence relation of isomorphism, considered as a subset of $\Mod(T)^2$, is analytic.
Given an $\L$-theory $T$ and an $\L'$-theory $S$, Friedman and Stanley say {\em $\Mod(T)$ is Borel reducible to $Mod(S)$,}
written $T\le_B S$, if there is a Borel function $f:\Mod(T)\rightarrow \Mod(S)$ such that for all models $M,N\models T$, $f(M)\cong_{\L'} f(N)$ if and only if $M\cong_\L N$.
Say the theories $T,S$ are {\em Borel equivalent} if both $T\le_B S$ and $S\le_B T$.  Thus, we get a notion of {\em Borel cardinality} comparing different classes of isomorphism types of countable structures.
Friedman and Stanley proved that there is a maximal class of theories with respect to Borel reducibility.  Theories in this maximal class are called {\em Borel complete}.

With this terminology, we seek to determine the Borel complexity of $Mod(T_R)$, the class of left $R$-modules with universe 
$\omega$.  Here, $T_R$ is the set of first order sentences in the language $\L_R=\{+,0,\cdot_r\}_{r\in R}$ axiomatizing `being a 
left $R$-module.'    The theory $T_R$ is incomplete, but there are natural completions.  As notation, 

$$R^{(\omega)}=\bigoplus_{\omega} \RsubR \quad\hbox{and} \quad R^\omega=\prod_\omega \RsubR$$
where $\RsubR$ is the left $R$-submodule with universe $R$,  $R^{(\omega)}$ is the (countable) direct sum, while the full product $R^\omega$ has size continuum.  
The modules $R^{(\omega)}$ and $R^\omega$  are  elementarily equivalent and both satisfy Fact~\ref{product}.

In \cite{LU}, the authors gave a complete solution to the Borel completeness of $T_R$ when the countable ring $R$ is commutative.

\begin{Theorem}[Theorem 1.2 of \cite{LU}]  Let $R$ be any countable, commutative ring (with unit). Then the following are equivalent:
\begin{enumerate}
		\item $R$ is an Artinian principal ideal ring;
		\item $R$ is of finite representation type (every module is a direct sum of indecomposable modules, and there are only finitely many indecomposable modules up to isomorphism);
		\item For every $R$-module $M$, $\mbox{Th}(M)$ is $\omega$-stable;
		
				\item There are only countably many isomorphism classes of countable $R$-modules;
				\item  $T_R$ is not Borel complete.  
			\end{enumerate}
\end{Theorem}

However, in \cite{LU}, the case for non-commutative $R$ was largely unexplored.
The current paper is devoted to handling the case when the countable ring $R$ has an infinite,  strictly descending chain
of pp-definable subgroups.  
(The pp stands for `positive primitive' and is defined in Section~2.)  Classes of countable modules of such rings have been studied before.  
For example, Bass's Theorem P  \cite{Bass} shows that every non-left-perfect ring has a strictly descending chain of principal right ideals, and hence a strictly descending chain of pp-definable subgroups.  Also, 
Prest~\cite{Prest} proves that for any countable ring $R$, $Th(\RsubR)$ is not $\omega$-stable if and only if $R$ has such a chain of pp subgroups.  

\begin{Theorem}  \label{big2}  Let $R$ be any countable ring with unit.  If the associated left $R$-module $\RsubR$ admits
a strictly descending chain of pp-definable subgroups, then the complete theory $Th(R^{(\omega)})$ 
is Borel complete.   Thus, both $Th(R^\omega)$ and  the incomplete theory $T_R$ are Borel complete as well. 
\end{Theorem}

Theorem~\ref{big2} directly generalizes Theorems~3.6 and 3.24 from \cite{LU}, but its proof is much more involved because we do not assume commutativity of $R$.   The short explanation of how commutativity was used in \cite{LU} is that $G_\infty=\bigcap_{n\in\omega} G_n$ is a 2-sided ideal of $R$.  The proofs of Theorems~3.6 and 3.24 both begin by modding out by this ideal.     

As $G_\infty$ is not typically a 2-sided ideal, with Definition~\ref{L} we define a surrogate 2-sided ideal $L^R$ that depends not only on $R$, but on the sequence $(G_n:n\in\omega)$.   
We also need a way of canonically closing a countable subset of a module to a module in a highly controlled way.   
We are interested in types over sets whose `positive part' is given by a single pp formula.  
With Definition~\ref{fgdef} we call such types {\em finitely generated} (f.g.) and we
introduce the notion of an {\em f.g.\ hull of a set}.  Under the assumptions of $T=T^{\aleph_0}$ and countability, with Proposition~\ref{uniquehull}  we prove

\begin{Proposition}  Suppose $T=T^{\aleph_0}$, the ring $R$ is countable, and $A\subseteq M$ is any countable subset of a module $M\models T$.
Then an f.g.\ hull $N\supseteq A$ exists and is unique up to isomorphism over $A$.
\end{Proposition}
\noindent We expect that the notion of an f.g.\ hull will be of independent interest beyond this paper.  

Recall that a ring $R$ is {\em left perfect} if all left $R$-modules have a projective cover.  These were studied by Bass~\cite{Bass},
who proved many equivalents of this notion.  The following corollary is proved in Section~\ref{last}.

\begin{Corollary}  \label{BassCor}   If a countable ring $R$ is not left perfect, then $T_R$ is Borel complete.
\end{Corollary}

However, Corollary~\ref{BassCor} is not a characterization.   There do exist countable, left perfect rings $R$ for which $T_R$ is Borel complete.\footnote{One commutative example is the finite ring $R={\mathbb F}_2[x,y]/(x^2,y^2,xy)$.  Being finite, $R$ is left perfect.   But  the ideals $(x)\cap(y)=(0)$ and $Ann(x)=Ann(y)=\{ax+by:a,b\in {\mathbb F}_2\}$, so by Theorem~3.21 of \cite{LU}, $T_R$ is Borel complete.}
However, among simple rings $R$, our results here do give a complete characterization, answering a question from \cite{LU}.  
Recall that a ring $R$ is {\em simple}  if it has no proper 2-sided ideals.    This corollary is also proved in Section~\ref{last}.

\begin{Corollary}  \label{simplecase}
Suppose $R$ is a countable simple ring with unit.  The following are equivalent.
\begin{enumerate}
\item  $T_R$ is not Borel complete;
\item  $R$ is simple artinian;
\item  $R\cong M_n(D)$ for some division ring $D$.
\end{enumerate}
\end{Corollary}

A final corollary of Theorem~\ref{big2} is that it sharpens existing proofs of TFAB, the theory of torsion-free abelian groups
given in \cite{PS} or \cite{LU}.   In particular, we now see that the complete theory $Th(\Z^{(\omega)})$ (all of whose models are
in TFAB) is Borel complete.

\subsection{Overview of the arguments}
Akin to the strategy of proving NP-completeness of some combinatorial statement, in order to prove the Borel completeness of some $\L$-theory $T$, it suffices to find a Borel reduction from some theory $S$, possibly in a different language, that is known to be Borel complete.

Here, we quote another main result from \cite{LU}.  For $R$ any ring, a {\em tagged left $R$-module} $(M,V_n:n\in\omega)$
is a left $R$-module $M$, and countably many left $R$-submodules $V_n$.  This notion was introduced by G\"obel and Shelah in \cite{GS}.
In the expanded language 
$\L'_R=\{+,0,\cdot_r,U_n:r\in R,n\in\omega\}$ we consider the theory $T_R'$ of all tagged left $R$-modules.  

\begin{Theorem}[Theorem~1.4 of \cite{LU}]   For any ring $R$ (not necessarily commutative) the theory $T_R'$ of tagged left
$R$-modules is Borel complete.
\end{Theorem}

The na\"ive approach toward proving Theorem~\ref{big2} would be to find a Borel reduction from $T_R'$ to $T_R$, whenever
$R$ has a distinguished descending chain $(G_n:n\in\omega)$ of pp-definable subgroups. 
Following the template of Theorems~3.6 and 3.24 of \cite{LU}, we form the inverse limit $\F=\varprojlim R/G_n$
and find a `sufficiently generic' countable sequence $\gammabar=(\gamma_k:k\in\omega)$.   With this sequence in hand, we would use the sequence $\gammabar$ to encode the distinguished subspaces of a given tagged $R$-module $(M,V_n)$ into a single $R$-module.    In order to show that this mapping preserves isomorphisms in both directions, we would like to use omitting types.  
However, for this approach, there is  a certain amount of information loss.
For the types $p_\gamma(x,a)$ we use, there is a uniformly defined $R$-submodule
$L^N$ of any left $R$-module $N$ (independent of $\gamma$) so that if $p_\gamma(x,a)$ is realized in $N$, then so is $p_\gamma(x,a+\ell)$ for every $\ell\in L^N$ (see Lemma~\ref{typesmodL}).   The good news is that $L^R$, the corresponding submodule of $\RsubR$, is actually a proper, 2-sided ideal of $R$ (Lemma~\ref{2sided}).  Thus, when we finally exhibit a Borel reduction, it will be from tagged $R/L^R$-modules to elementary submodels of a carefully constructed countable $R$-module $N^*$.


In addition to proving the theorem, there are two novel module-theoretic notions introduced here.   The first is the proper 2-sided ideal $L^R$ mentioned above, which depends not only on the ring $R$ but also on the specific choice of a descending sequence 
$(G_n:n\in\omega)$.   As the sequence will be fixed throughout, we do not explicitly mention the sequence in the notation 
$L^R$.  

The second, arguably more important contribution is the notion of f.g.\ hulls over a set.   Whenever $R$ is a countable ring and $T$ is a theory of left $R$-modules satisfying 
$T=T^{\aleph_0}$\footnote{The theory of the direct sum $R^{(\omega)}$ satisfies $T=T^{\aleph_0}$ for any choice of $R$.}
then for any countable subset $A\subseteq M$ of any module $M\models T$, there is a countable f.g.\ hull $N\supseteq A$
that is unique up to isomorphism over $A$ (see Proposition~\ref{uniquehull}).  This notion acts as a surrogate for a `prime model over $A$' which need not exist in a non-$\omega$-stable setting.  
$N$ being an f.g.\ hull over $A$ is a conjunction of two conditions.  Modules $M\supseteq A$ satisfying the first condition, that  every finite tuple $\mbar$ has $\tp(\mbar/A)$ finitely generated over $A$ were studied by Kucera and Rothmaler~\cite{KR}.  They showed that any such $M$ is pure-projective over $A$ and 
that $M$ was `constructible' over $A$.  However, that condition alone is not sufficient to imply uniqueness over $A$.  The uniqueness is guaranteed by adding a second `saturation' clause, stating that for every finite tuple $\nbar$ from $N\supseteq A$, every f.g.\ 1-type over $A\nbar$ realized in the monster model is realized in $N$.

Once the Preliminaries and the discussion of f.g.\ hulls are dispensed with, in Section~\ref{fixring} we fix, once and for all, a countable ring $R$ with a fixed, strictly descending sequence $(G_n:n\in\omega)$ of pp-definable subgroups.  We form the inverse limit $\F=\varprojlim
R/G_n$ and define and develop both the 2-sided ideal $L^R$ and
the family of pp-types $\{p_\gamma(x,y):\gamma\in\F\}$.    This culminates in Lemma~\ref{typesmodL}, which places a constraint on how independent the types $p_\gamma(x,y)$ are.   

In Section~\ref{transfer}, we develop general machinery for producing a `suitably generic sequence' $\gammabar$ for a given countable, tagged $R$-module $M$.     However, in what follows, we use these results only once.   We construct a single generic sequence $\gammabar$ that works with the  universal free-like tagged $R$-module $F^*$.     
In Section~\ref{modLRX} we prove Proposition~\ref{almost}, which is  an `approximation' to a Borel reduction from tagged $R$-modules to $R$-modules.
  We use the uniqueness of f.g.\ hulls to show that isomorphic tagged $R$-modules with universe $RX$ have isomorphic images, but the reverse direction, which relies on modules realizing or omitting certain types $p_\gamma(x,a)$
   is only ``modulo $L^{RX}$."  
  
In Section~\ref{Borelreduction} we prove that, letting $S=R/L^R$, there is a true Borel reduction from tagged $S$-modules to 
countable models of $Th(R^{(\omega)})$.

Finally, in Section~\ref{last} we prove the two Corollaries mentioned above.

\section{Preliminaries}

Fix any ring $R$ (with unit) and we study the class of left $R$-modules.  Model-theoretically, these modules form  an elementary class in the language $L_R=\{+,-,\cdot_r\}_{r\in R}$.  Here, we skim the main model theoretic results, but the reader should consult either Prest \cite{Prest} or Ziegler \cite{Ziegler} for proofs and more details.  

The most important $L_R$-formulas are the {\em positive primitive (pp)} ones
$$\phi(x_1,\dots,x_n):=\exists y_1\dots \exists y_k\bigwedge_{\ell=1}^p\left( \sum_{i=1}^n r_{i,\ell} x_i +\sum_{j=1}^k s_{j,\ell} y_j=0\right)$$
For a fixed $R$-module $M$, $\phi(M)=\{(a_1,\dots,a_n)\in M^n:M\models \phi(a_1,\dots,a_n)\}$ is called a {\em pp-definable subgroup}
and it asserts the tuples $\abar\in M^n$ for which the finite system of homogeneous $R$-linear equations has a solution in $M$.

The lattice of pp-definable subgroups $\phi(M)$ is highly relevant.  The elementary theory of $M$, $Th(M)$, is determined by the indices $[\psi(M):\phi(M)\cap\psi(M)]\in\{1,2,3,\dots,\infty\}$ as $\phi(\xbar)$, $\psi(\xbar)$ range over all pp formulas.

Of particular interest to us will be the modules $N$ for which each of these indices are either 1 or $\infty$.  Clearly, for any 
$R$-module $M$, either the infinite direct sum $M^{(\omega)}$ or the product $M^\omega$ have this additional property.   Following Garavaglia \cite{Gar} and Pillay-Prest~\cite{PP}, we say $T=T^{\aleph_0}$ if $T$ is a theory of modules for which this property holds.   For us, we use that  the theory $T=Th(R^{(\omega)})$ has this property.  

The major result is the Baur-Monk quantifier elimination for theories of modules.  For any $R$-module $M$, every $L_R$-formula
$\theta(\xbar)$ is equivalent to a boolean combination of pp formulas in the variables $\xbar$. 

Frequently, we add parameters to supplant some free variables.  That is, for an $R$-module $M$, $\theta(\xbar,\zbar)$ any formula, and $\bbar\in M^{\lg(\zbar)}$, then $\theta(\xbar,\bbar)$ is a formula with parameters.  If $\phi(\xbar,\zbar)$ is pp, then
the solution set $\phi(M,\bbar)$ is either empty, or else is a coset of the pp-definable subgroup $\phi(M,\0bar)$, formed by substituting $0's$ for every $z\in\zbar$.  

With this notion, if $M$ is a submodule of $N$, then $M$ is {\em pure in $N$} if, for every finite $\abar$ from $M$, every pp-formula $\phi(\xbar,\abar)$ with a solution in $N$ has a solution in $M$.  It follows from the Baur-Monk quantifier elimination that
{\em $M$ is an elementary substructure of $N$,} written $M\preceq N$ if they are elementarily equivalent, written $M\equiv N$, and $M$ is pure in $N$.

Once we settle on a theory $T$ of modules, we fix a {\em monster model $\C\models T$,} which is very large and very saturated.  
Whenever we talk about a model, we will mean an elementary substructure $M\preceq \C$ of the monster.  

Suppose $A\subseteq\C$ is a small subset (relative to the saturation of $\C$).  
A {\em pp partial type $\Gamma^+(\xbar,A)$} is a set of pp formulas $\phi(\xbar,\abar)$, each with free variables $\xbar$ and
with $\abar$ from $A$, that is {\em satisfiable}, i.e., there is some $\bbar\in\C^{\lg(\xbar)}$ {\em realizing $\Gamma^+(\xbar,A)$,}
i.e., $$\bbar\in\bigcap\{\phi(\C,\abar):\phi(\xbar,\abar)\in \Gamma^+(\xbar,A)\}$$
Checking the satisfiability of a pp partial type is expedited by the Compactness Theorem.  
$\Gamma^+(\xbar,A)$ is satisfiable if and only if every  pp formula $\phi(\xbar,\abar)$
formed as finite conjunction of elements of $\Gamma^+(\xbar,A)$ has a solution in $\C$.  

\subsection{Consequences of $T=T^{\aleph_0}$}

Here, let $R$ be any ring and $\C\models T$ be a monster model of any  theory of left $R$-modules satisfying $T=T^{\aleph_0}$, i.e., 
for all pp formulas $\phi(\xbar),\psi(\xbar)$, 
either $\psi(\C)\le\phi(\C)$ or the index 
$[\psi(\C):
   \varphi(\C)\cap\psi(\C)]=\infty$.

The following weak version of Neumann's Lemma will be helpful.

\begin{Fact}[Weak Neumann's Lemma]  \label{Neumann} 
Suppose $C$ is a coset of a group $G$ that is covered by finitely many cosets $D_i$ of subgroups $H_i\le G$.   Then some $H_i$ has finite index in $G$.
\end{Fact}

\begin{Definition}  {\em  Assume $A\subseteq B$ and $\Gamma^+(\xbar,A)$ and $\Delta^+(\xbar,B)$ are partial pp-types with the same free variables $\xbar$.  We say {\em $\Delta^+(\xbar,B)$ extends
 $\Gamma^+(\xbar,A)$}  if $\Gamma^+(\xbar,A)\subseteq\Delta^+(\xbar,B)$.  Clearly, in this case, $\Delta^+(\C,B)\subseteq \Gamma^+(\C,A)$.  In particular, $\Delta^+(\C,\0bar)\le \Gamma^+(\C,\0bar)$ as subgroups.  
 
Conversely, we say {\em $\Gamma^+(\xbar,A)$ entails $\Delta^+(\xbar,B)$,} written $\Gamma^+(\xbar,A) \models\Delta^+(\xbar,B)$, if $\Gamma^+(\C,A)\subseteq \Delta^+(\C,B)$.
By compactness, $\Gamma^+(\xbar,A)\models \Delta^+(\xbar,B)$ if and only if, for every $\delta(\xbar,\bbar)\in \Delta^+(\xbar,B)$ there is a finite conjunction $\theta(\xbar,\abar)$ from
$\Gamma^+(\xbar,A)$ such that $\theta(\C,\abar)\subseteq \delta(\C,\bbar)$.

For any set $B$, a {\em complete type over $B$,}  $p(\xbar)$, is a satisfiable set of $L(B)$ formulas, all of whose free variables are from $\xbar$ such that, for every $L(B)$-formula $\psi(\xbar,\bbar)$, exactly one of $\psi(\xbar,\bbar),\neg\psi(\xbar,\bbar)\in p$.
The Stone space $S(B)$ consists of all complete types over $B$.  For every $p(\xbar)\in S(B)$,
$p^+(\xbar)=\{\hbox{all pp formulas $\phi(\xbar,\bbar)\in p$}\}$.

Suppose $A\subseteq B$, $\Gamma^+(\xbar,A)$ is a partial pp-type, and $p\in S(B)$ is a complete type.  We say {\em  $\Gamma^+(\xbar,A)$ generates $p$} if
$\Gamma^+(\xbar,A)\subseteq p^+$ and entails $p^+$.  
}
\end{Definition}

\begin{Proposition}   \label{main}    Suppose $T=T^{\aleph_0}$.   For every partial pp-type $\Gamma^+(\xbar,A)$ and every $B\supseteq A$,
there is a unique $p\in S(B)$ generated by $\Gamma^+(\xbar,A)$.   
\end{Proposition}

\begin{proof}  
Given $A\subseteq B$ and partial pp-type $\Gamma^+(\xbar,A)$, let
$$\Phi(\xbar,B):=\{\neg\delta(\xbar,\bbar):\hbox{$\delta(\xbar,\bbar)$ pp and} \  \Gamma^+(\xbar,A)\not\models \delta(\xbar,\bbar)\}$$

\medskip
\noindent{\bf Claim 1.}  $\Gamma^+(\xbar,A)\cup\Phi(\xbar,B)$ is finitely satisfiable.

\begin{proof}  Without loss, $\Gamma^+(\xbar,A)$ is closed under finite conjunctions, so choose $\theta(\xbar,\abar)\in\Gamma^+(\xbar,A)$ and choose $\{\neg\delta_i(\xbar,\bbar_i):i<n\}\subseteq \Phi(\xbar,B)$.  We may assume that $\theta(\xbar,\abar)\wedge\delta_i(\xbar,\bbar)$ is satisfiable for each $i<n$, since otherwise we could remove it from the conjunct.
Also, by replacing $\neg\delta_i(\xbar,\bbar_i)$ by $\neg(\theta(\xbar,\abar)\wedge\delta_i(\xbar,\bbar_i))$, we may assume $\delta_i(\C,\bbar_i)\subseteq\theta(\C,\abar)$.  
If equality held for any $i$, we would have $\theta(\xbar,\abar)\models\delta_i(\xbar,\bbar_i)$, contradicting $\delta_i(\xbar,\bbar_i)\in \Phi(\xbar,B)$.  Thus, since $T=T^{\aleph_0}$, Neumann's Lemma~\ref{Neumann} implies that $\bigcup_{i<n}\delta_i(\C,\bbar_i)\neq \theta(\C,\abar)$.  Thus, 
$\C\models\exists\xbar(\theta(\xbar,\abar)\wedge\bigwedge_{i<n} \neg\delta_i(\xbar,\bbar_i))$.
\qed\end{proof}

Now choose $\ebar\in\C$ realizing $\Gamma^+(\xbar,A)\cup\Phi(\xbar,B)$ and let $p(\xbar):=\tp(\ebar/B)\in S(B)$.  To see that $\Gamma^+(\xbar,A)$ generates $p$, it suffices to show that every pp $\delta(\xbar,\bbar)\in p^+$ is entailed by $\Gamma^+(\xbar,A)$.   Choose any $\delta(\xbar,\bbar)\in p^+$.  If it were not entailed by $\Gamma^+(\xbar,A)$, then we would
have $\neg\delta(\xbar,\bbar)\in \Phi(\xbar,B)\subseteq p$, which is a contradiction.

Finally, to prove uniqueness, suppose  $p'\in S(B)$  is also generated by $\Gamma^+$.   Then, for
every pp formula 
$\delta(\xbar,\bbar)$,
$$
 \delta(\xbar,\bbar)\in p'
 \quad\Longleftrightarrow 
 \quad
 \Gamma^+\models\delta(\xbar,\bbar)
 \quad\Longleftrightarrow 
  \delta(\xbar,\bbar)\in p$$
That is, $(p')^+=p^+$.   By the Baur-Monk quantifier elimination we conclude $p'=p$.
\qed
\end{proof}

\begin{Remark}  {\em In fact, the Proposition characterizes theories of modules satisfying $T=T^{\aleph_0}$.   Indeed, suppose there were pp-definable groups $\phi(\C)\le\psi(\C)$ with $1<[\psi(\C):\phi(\C)]<\aleph_0$.   Put $A:=\emptyset$ and take $B$ to be any elementary substructure $M\preceq \C$.   Choose $S=\{b_1,\dots,b_n\}\subseteq B$ to be a set of representatives for the cosets of $\phi(\C)$ inside $\psi(\C)$.  Note that any complete type $p(x)\in S(M)$ containing $\psi(x)$
must contain some pp-formula $\phi(x-b_i)$ for some $b_i\in S$.   Thus, the partial type $\Gamma(x):=\{\psi(x)\}$ does not generate any complete type $p(x)\in S(M)$.
}
\end{Remark}

\begin{Lemma}  \label{typegenerates}   Suppose $T=T^{\aleph_0}$.   
\begin{enumerate}
\item  If $\theta(\ybar,\abar)$ and $\{\psi_j(\ybar,\bbar_j):j<m\}$ are pp formulas
and $\theta(\ybar,\abar)\models\bigvee_{j<m} \psi_j(\ybar,\bbar_j)$, then
$\theta(\ybar,\abar)\models\psi_j(\ybar,\bbar_j)$ for at least one $j<m$.
\item  If a complete type $p(\ybar)\in S(A)$ entails $\bigvee_{j<m}\gamma_j(\ybar,\bbar_j)$, where each $\gamma_j$ is pp, 
 then $p^+(\ybar)\models\gamma_j(\ybar,\bbar_j)$ for some $j<m$.  
\end{enumerate}
\end{Lemma}

\begin{proof}  
(1)  If $\theta(\ybar,\abar)$ is not satisfiable, there is nothing to prove, so assume $\theta(\C,\abar)\neq\emptyset$.  Similarly, we may assume each $\psi_j(\ybar,\bbar_j)$ is satisfiable as well.  
By replacing each $\psi_j$ by $\theta\wedge\psi_j$, we may also assume $\psi_j(\C,\0bar)\le\theta(\C,\0bar)$ for each $j<m$.   If $\theta(\ybar,\abar)$ did not entail any $\psi_j(\ybar,\bbar_j)$, then by $T=T^{\aleph_0}$, we would have the index $[\theta(\C,\0bar):\psi_j(\C,\0bar)]$ infinite for every $j<m$.  But this contradicts Neumann's Lemma.

(2)  By compactness, choose $\eta(\ybar,\abar)\in p$ such that $\eta(\ybar,\abar)\models\bigvee_{j<m}\gamma_j(\ybar,\bbar_j)$.  
By Baur--Monk elimination, choose a formula $\chi(\ybar)\in p$ of the form
$$
        \chi(\ybar)
        =
        \theta(\ybar,\abar)
        \wedge
        \bigwedge_{i<n}\neg\psi_i(\ybar,\abar_i)$$
such that \(\chi\models\eta\) and  \(\theta\) and the \(\psi_i\)'s
are pp.  Then
$
 \theta(\ybar,\abar)
 \models
 \bigvee_{j<m}\gamma_j(\ybar,\bbar_j)
 \ \vee\
 \bigvee_{i<n}\psi_i(\ybar,\abar_i).$
By (1), $\theta$ entails one of the displayed disjuncts.
It cannot entail any $\psi_i$, since $\chi$ is consistent.
Hence $\theta$ entails some
$\gamma_j(\ybar,\bbar_j)$, as required.
\qed
\end{proof}

For any stable theory $T$, if $A\subseteq B$, then the consistent $L(B)$-formula $\phi(\xbar,\bbar)$  forks over $A$ if and only if $\phi(\xbar,\bbar)$ {\em divides over $A$}, i.e., there is an infinite set $\{\bbar_i:i\in\omega\}$ and an integer $k\ge 2$ such that $\tp(\bbar_i/A)=\tp(\bbar/A)$ for every $i$ and the set $\{\phi(\C,\bbar_i):i\in\omega\}$ is $k$-inconsistent.   For theories of modules we can say more.  Since conjugate cosets are disjoint or equal, a  consistent pp-formula $\phi(\xbar,\bbar)$ divides over $A$ if and only if the orbit of $\phi(\C,\bbar)$ is infinite under $Aut(\C/A)$, the group of automorphisms of $\C$ fixing $A$ pointwise.   For theories of modules with $T=T^{\aleph_0}$ we have the following.  In different terminology, this was also proved by Pillay and Prest~\cite{PP}.


\begin{Corollary}  \label{nonforking}
Suppose $T=T^{\aleph_0}$, $A\subseteq B$ and $p\in S(A)$.   Then the type $q\in S(B)$ generated by $p^+$ is the unique non-forking extension of $p$ to $S(B)$.
\end{Corollary}

\begin{proof}  
Since theories of modules are stable, choose a nonforking extension
$r\in S(B)$ of $p$.  We show that $p^+$ generates $r$.  By
Proposition~\ref{main}, this will imply $r=q$.  Applying the same argument
to any nonforking extension will also give uniqueness.

Choose any pp formula $\delta(\xbar,\bbar)\in r$ and let $t(\ybar)=\tp(\bbar/A)$.
Let $\gamma(\ybar,\bbar)$ be the pp formula $\exists\xbar (\delta(\xbar,\ybar)\wedge\delta(\xbar,\bbar))$.   Note that for any $\bbar'$, $\gamma(\bbar',\bbar)$ holding implies that the cosets
$\delta(\C,\bbar')$ and $\delta(\C,\bbar)$ are equal.   

\medskip
\noindent{\bf Claim.}  Some $\theta(\ybar,\abar)\in t^+$ satisfies $\theta(\ybar,\abar)\models\gamma(\ybar,\bbar)$.

\begin{proof}  Since $\delta(\xbar,\bbar)$ does not divide over $A$, there are only finitely many distinct sets $\delta(\C,\bbar')$ for any $\bbar'$ realizing $t(\ybar)$.  Choose representatives
$\{\bbar_i:i<m\}\subseteq t(\C)$ with $\bbar_0=\bbar$ so that 
$$t(\ybar)\models\bigvee_{j<m} \exists \xbar (\delta(\xbar,\ybar)\wedge\delta(\xbar,\bbar_j))$$
As the consequent is a disjunction of  pp formulas, Lemma~\ref{typegenerates}(2) gives some $j^*<m$ such that 
$$t^+(\ybar)\models\exists \xbar(\delta(\xbar,\ybar)\wedge\delta(\xbar,\bbar_{j^*}))$$
As $\bbar$ realizes $t^+(\ybar)$, $\delta(\C,\bbar_{j^*})\cap\delta(\C,\bbar)\neq\emptyset$, hence $j^*=0$.  That is,  $t^+(\ybar)\models\gamma(\ybar,\bbar)$.
By compactness, there is a pp formula $\theta(\ybar,\abar)\in t^+$ entailing $\gamma(\ybar,\bbar)$.
%
%
\qed
\end{proof}

Continuing, let $\rho(\xbar,\abar):=\exists\ybar(\theta(\ybar,\abar)\wedge\delta(\xbar,\ybar))$.
Visibly, $\rho(\xbar,\abar)$ is pp, hence $\rho(\xbar,\abar)\in p^+$ by taking $\bbar$ to witness the existential quantifier.  It remains to show that $\rho(\xbar,\abar)\models \delta(\xbar,\bbar)$.
For this, choose any $\cbar'$ realizing $\rho(\xbar,\abar)$.   Choose $\bbar'$ such that $\theta(\bbar', \abar)$ and $\delta(\cbar',\bbar')$ hold.   By the Claim,  $\gamma(\bbar',\bbar)$ holds,
hence $\delta(\C,\bbar')=\delta(\C,\bbar)$, so $\delta(\cbar',\bbar)$ holds, as required.
\qed
\end{proof}


Given this correspondence, and the fact that any theory of modules is stable, we can adopt the notation of non-forking.

\begin{Notation} {\em  For any sets $A,B,C$, we write $\fg A B C$ if, for every finite $\abar$ from $A$, $\tp^+(\abar/B)$ generates
$\tp(\abar/BC)$.  
}
\end{Notation}

Given our correspondence, the  following Fact is just rewriting the fact that non-forking 
 satisfies Symmetry and Transitivity in any stable theory,  
\begin{Fact}   \label{useful}    Suppose $T=T^{\aleph_0}$.  
\begin{enumerate}
\item{\bf Transitivity} Suppose $B\subseteq C\subseteq D$.  Then for any finite tuple $\abar$,
$\tp^+(\abar/B)$ generates $\tp(\abar/D)$ if and only if $\tp^+(\abar/B)$ generates $\tp(\abar/C)$ and $\tp^+(\abar/C)$ generates $\tp(\abar/D)$.
\item{\bf Symmetry} For any $\abar,B,\cbar$,  $\tp^+(\abar/B)$ generates $\tp(\abar/B\cbar)$ if and only if $\tp^+(\cbar/B)$ generates 
$\tp(\cbar/B\abar)$.
\end{enumerate}
\end{Fact}  

%
%

\subsection{Model theory of free modules}

Suppose $R$ is a ring and $X$ is any  countably infinite set.   We write $RX$ for the {\em free $R$-module with basis $X$}, i.e.
$$RX=\{\sum r_ix_i:r_i\in R, x_i\in X\}$$
Any such $RX$ is isomorphic to $R^{(\omega)}$ and is elementarily equivalent to the full product $R^X$, which is isomorphic to
the full product $\prod_\omega \RsubR$.

The following Fact is just unpacking the definitions.  The same fact holds for $R^X$, but in this paper only use the countable modules $RX$.  

\begin{Fact}  \label{product}  Suppose $\abar=(a_1,\dots,a_n)\in RX^n$ with each $a_i=\sum r_{ij} x_j$.
Then for any pp formula $\phi(\xbar)$
$$RX\models \phi(a_1,\dots,a_n)\quad \hbox{if and only if} \quad \hbox{for every $x_j\in X$, $\RsubR\models 
\phi(\sum r_{1j},\sum r_{2j},\dots, \sum r_{nj}$)}$$
\end{Fact}  

For lack of a better place, we include the following fact.

\begin{Fact}  \label{surjective}   
Let  $f:M\rightarrow N$ be any left  homomorphism among modules.  For any pp formula $\phi(\xbar)$ and any
$\abar$ from $M$, if $M\models \phi(\abar)$, then $N\models\phi(f(\abar))$.  
\end{Fact}  

Putting these two together yields the following fact.

\begin{Fact}  \label{1suffices}  Let $\phi(y)$ be any pp-formula with no hidden parameters and suppose
$\RsubR\models \phi(1)$.  Then $\C\models\forall y\phi(y)$.
\end{Fact}  

\begin{proof}   Assume $\RsubR\models \phi(1)$.  Since $\phi(R)$ is a right $R$-module, it follows that
$\RsubR\models\phi(a)$ for every $a\in R$.  
Thus, by Fact~\ref{product}, $RX\models\forall y\phi(y)$ for any infinite $X$.
  Since $RX\preceq \C$ we finish.
\qed
\end{proof}

\section{Introducing f.g.\ hulls}

Throughout this section we assume $\C\models T$ is a monster model for some theory of modules satisfying
$T=T^{\aleph_0}$.

\begin{Definition}  \label{fgdef} {\em  A   complete type $p(\xbar)\in S(A)$ is {\em finitely generated (f.g.) over $A$}  if there is some pp formula $\theta(\xbar,\abar)\in p$ such that
$\{\theta(\xbar,\abar)\}$ generates $p$ (hereafter, we say $\theta(\xbar,\abar)$ generates $p$ in place of $\{\theta(\xbar,\abar)\}$).    \\ Call a set $B\supseteq A$ {\em f.g.\ atomic over $A$}
if $\tp(\bbar/A)$ is f.g.\ over $A$ for every finite $\bbar$ from $B$.
\\  $B$ is called {\em f.g.\ full over $A$} if, for every finite $\bbar$ from $B$, every f.g.\ type $p\in S(A\bbar)$ is realized in $B$.
\\  $B$ is an {\em f.g.\ hull over $A$} if it is both f.g.\ atomic and f.g.\ full over $A$.
}
\end{Definition}

The remarkable fact is that under the assumption of $T=T^{\aleph_0}$, the f.g.\ types are always dense.  That is:

\begin{Lemma}  \label{dense}  For any set $A$ and any consistent  $L(A)$-formula $\psi(\xbar,\abar)$, there is an f.g.\ type $p\in S(A)$ containing $\psi(\xbar,\abar)$.
\end{Lemma}

\begin{proof}  In light of the pp elimination of quantifiers, it suffices to show this for any consistent $L(A)$-formula $\psi(\xbar,\abar)$ of the form
$$\theta(\xbar,\abar)\wedge\bigwedge_{i<n} \neg \gamma_i(\xbar,\abar)$$
where $\theta$ and each $\gamma_i$ are pp formulas.  By Proposition~\ref{main}, let $p(\xbar)\in S(A)$ be the unique complete type generated by $\theta(\xbar,\abar)$.
The type $p(\xbar)$ is visibly finitely generated over $A$, so it suffices to show $\theta(\xbar,\abar)\wedge\neg \gamma_i(\xbar,\abar)\in p$ for every $i<n$.
However, if this did not hold, then for some $i$, the pp-formula $\theta(\xbar,\abar)\wedge\gamma_i(\xbar,\abar)$ would be an element of $p^+$.  But, since $\theta(\xbar,\abar)$
generates $p$, we would have $\theta(\xbar,\abar)\models \gamma_i(\xbar,\abar)$, contradicting the fact that $\psi(\xbar,\abar)$ is consistent.
\qed
\end{proof}

 \begin{Lemma}  \label{TV}   For any set $A$, if $B\supseteq A$ is an f.g.\ hull over $A$, then $B$ is an elementary substructure of $\C$.
 \end{Lemma} 
 
 \begin{proof}  This is immediate by  the Tarski-Vaught test.  Take any consistent $L(B)$-formula $\psi(\xbar,\bbar)$ and by Lemma~\ref{dense}, choose a complete type
 $p\in S(A\bbar)$ that is f.g.\ over $A\bbar$ containing $\psi(\xbar,\bbar)$.   Since $B$ is f.g.\ full, there is some $b^*\in B$ realizing $p$.  Thus, $B$ contains a 
 realization of $\psi(\xbar,\bbar)$.
 \qed
 \end{proof}
 
 The following Lemma appears as Lemma~2.5  in \cite{KR}.  
 
 \begin{Lemma}  \label{pieces}  Suppose $A,\bbar,\cbar$ are given.   Then $\tp(\bbar\cbar/A)$ is f.g.\ if and only if both $\tp(\bbar/A\cbar)$ and $\tp(\cbar/A)$ are f.g.
 \end{Lemma}
 
 In light of Lemmas \ref{dense} and \ref{pieces}, it is natural to define an {\em f.g.-construction sequence over $A$} as an ordinal-indexed sequence $\<a_i:i<\alpha\>$ such that, 
 for every $i<\alpha$, $\tp(a_i/A_i)$ is finitely generated, where $A_i=A\cup\{a_j:j<i\}$.   As noted in \cite{KR}, an easy induction on $\alpha$ shows that every f.g.-construction sequence is f.g.-atomic over $A$.

 \begin{Proposition}  \label{hullexistence}   For any set $A$, there is an f.g.\ hull $M\supseteq A$ of cardinality at most $\kappa:=|R|+|A|+\aleph_0$.
 Moreover if $B\supseteq A$ is f.g.\ atomic over $A$ with $|B|\le\kappa$, then there is an f.g.\ hull $M$ over $A$ of size at most $\kappa$ with $M\supseteq B$.
 \end{Proposition}

 \begin{proof}      The first sentence obviously follows from the second, so fix $A\subseteq B$ with $B$ f.g.\ atomic over $A$.  To extend $B$ to an f.g.\ hull $M$ over $A$, all we need is a good secretary.   The key observation is that for any $B'\supseteq B$ of size at most $\kappa$, there are only $\kappa$ formulas $\theta(x,\bbar')$, hence there are only $\kappa$ complete f.g.\ types over $B'$.  
Thus, we can form an f.g.-construction sequence $\<b_i<\kappa\>$ over $B$ in such a way that for every finite $\bbar$ from 
 $B_\kappa:=B\cup\{b_i:i<\kappa\}$,  every f.g.\ type $p\in S(B\bbar)$ is realized in $B_\kappa$.    This is possible by repeated applications of Lemma~\ref{dense} and proper bookkeeping.   
 Thus,  the resulting $B_\kappa$ is f.g.\ full over $B$.   We argue that $B_\kappa$ is f.g.\ full over $A$ as well.   For this, choose any finitely generated $q\in S(A\bbar')$ for some $\bbar'$ from $B_\kappa$.  Let $q^*\in S(B\bbar')$ be the unique non-forking extension of $q$.  By the uniqueness of Proposition~\ref{main} and Corollary~\ref{nonforking}, $q^*$ is also finitely generated, hence realized in $B_\kappa$.
 
Furthermore, being the universe of an f.g.\ construction sequence over $B$, $B_\kappa$ is also f.g.\ atomic over $B$.  Since $B$ is f.g.\ atomic over $A$, it follows easily from Lemma~\ref{pieces} that  $B_\kappa$ is an f.g.\ hull over $A$.  That $B_\kappa\preceq\C$ is Lemma~\ref{TV}.
 \qed
 \end{proof}

 For the remainder of this section, in addition to $T=T^{\aleph_0}$, we additionally assume the ring $R$ is countable.

 \begin{Lemma}   \label{Borelpiece}   Suppose $N^*$ is a countable f.g.\ hull over $B$.  Then, for any $A\subseteq B$, there is $N_A\preceq N^*$ with $N_A$ an f.g.\ hull over $A$ and
 $\fg {N_A} A B$.  Moreover, the mapping $f:\P(B)\rightarrow \P(N^*)$ given by $A\mapsto N_A$ can be chosen to be  Borel.  
 \end{Lemma}
 
 \begin{proof}   Fix, for the whole of the argument, an enumeration $(d_\ell:\ell\in \omega)$ of $N^*$.  
 Fix any $A\subseteq B$.   We construct $N_A$ as an f.g.\ construction sequence in $\omega$ steps as in the proof of Proposition~\ref{hullexistence}.   Suppose that at stage $i<\omega$,
 we have constructed $A_i=A\cup\<a_j:j<i\>$ and $\fg {A_i} A B$ holds.
 Now, suppose we are given a finitely generated type  $p=\tp(a_i/A_i)$ and we want to add a realization of $p$ in $N^*$ to $A_{i+1}$ that will maintain the independence.  Moreover, this should be done in a canonical way to obtain the Borelness of the map $f$.  
   Suppose $p$ is generated by $\theta(x,\abar)$ and let  $p'\in S(A_i\cup B)$ be the unique type in $S(A_i\cup B)$ generated by
   $\theta(x,\abar)$.      
   As 
 $N^*$ is f.g.\ full,  and as $A_i\setminus B$ is finite,  $p'$ is realized in $N^*$.
Put $a_{i}:=d_\ell$, where $\ell$ is least such that $d_\ell$ realizes $p'$ and let $A_{i+1}=A_ia_i$.
 Since $p'$ is the non-forking extension of $p$ we have $\fg {a_i} {A_i} B$.
 By symmetry and transitivity we obtain $\fg {A_{i+1}} A B$, and we continue.  With good bookkeeping, the resulting
 $A_\omega$ will be a suitable $N_A$.  
 \qed
 \end{proof}
 
 Now, if both the ring $R$ and the set $A$ are countable, we obtain a canonical object.  
 
 \begin{Proposition}   \label{uniquehull}  Suppose $T=T^{\aleph_0}$ and both $R$ and $A$ are countable.  There exists an f.g.\ hull $N^*$ that is unique up to isomorphism over $A$.  
 \end{Proposition}
 
 \begin{proof}  The existence of a countable f.g.\ hull  $N^*$ follows immediately from Proposition~\ref{hullexistence}.  Now suppose $N_1,N_2$ are both countable f.g.\ hulls over $A$.
 Let $(\F,\le)$ consist of all partial  $L(A)$-elementary functions $f:\bbar\mapsto\cbar$ with $\bbar\subseteq N_1$, $\cbar\subseteq N_2$, ordered by inclusion.  
 We claim that $(\F,\le)$ is a back and forth system, which suffices as both $N_1,N_2$ are countable.  By symmetry, it suffices show `forth.'  So, choose $d\in N_1$ arbitrarily.
 As $N_1$ is f.g.\ atomic,   $\tp(\bbar d/A)$ is finitely generated.  By Lemma~\ref{pieces}, both $p:=\tp(\bbar/A)$ and $q:=\tp(d/A\bbar)$ are f.g.\ as well.  As $f$ is $L(A)$-elementary,
 the type $f(q)\in S(A\cbar)$ is f.g.\ over $A\cbar$.  Since $N_2$ is f.g.\ full, there is a realization $e\in N_2$ of $f(q)$.  Then the one-point extension $g:=f\cup\{(d,e)\}\in \F$ as required.
 \qed
 \end{proof}

 \begin{Remark}  {\em  Even in the countable case, the f.g.\ hull  $N^*$ over $A$ constructed above is not really a `hull' at all.  That is, there is a proper elementary substructure $N'\preceq N^*$ that is also an f.g.\ hull over $A$.  To see this, first note that by Proposition~\ref{uniquehull}, showing the existence of $N'$ is equivalent to showing the existence
 of a proper countable elementary extension  $N^{**}$ of $N^*$ that is also an f.g.\ hull over $A$.   To construct  such a $N^{**}$, first take $\theta(x)$ to be $x=x$
 and let $p\in S(N^*)$ be the complete type generated by $\theta(x)$.  Let $c$ be any realization of $p$.  As  $x\neq n\in p$ for every $n\in N^*$, $c\not\in N^*$.  But $N^* c$ is f.g.\ atomic over $A$, hence by the Moreover clause of Proposition~\ref{hullexistence}, there is a countable f.g.\ hull $N^{**}$ over $A$ containing $N^*c$, as required.
 }
 \end{Remark}

 \begin{Lemma}  \label{onestepfg}  Suppose $A\subseteq B$, $M$ is an f.g.\ hull over $A$, $N$ is an f.g.\ hull over $MB$,
 and $\fg M A B$.   Then $N$ is an f.g.\ hull over $B$.
 \end{Lemma}
 
 \begin{proof}  To show that $N$ is f.g.\ atomic over $B$, choose $\cbar$ from $N$.   Since $N$ is an f.g.\ hull over $MB$, choose
 a pp-formula $\delta(\xbar,\mbar,\bbar)$ generating $\tp(\cbar/MB)$.   Since $M$ is f.g.\ atomic over $A$, choose a pp-formula $\phi(\xbar,\abar)\in\tp^+(\mbar/A)$ generating $\tp(\mbar/A)$.  By the independence and Fact~\ref{useful}, $\phi(\xbar,\abar)$ generates $\tp(\mbar/B)$.   Thus, by Lemma~\ref{pieces},  $\tp(\cbar\mbar/B)$ is finitely generated, so 
 $\tp(\cbar/B)$ is also finitely generated by a second application of Lemma~\ref{pieces}.  
 
 Showing $N$ f.g.\ full over $B$ does not require independence.   Choose any $\dbar$ from $N$ and fix any
 finitely generated $p(x)\in S(B\dbar)$.  Say the pp-formula $\delta(x,\dbar,\bbar)$ generates $p$.  By Proposition~\ref{main}
 let $p^*\in S(MB\dbar)$ be generated by $\delta(x,\dbar,\bbar)$.   By uniqueness, it follows that $p^*\supseteq p$.
 As $N$ is f.g.\ full over $MB$, there is an $e\in N$ realizing $p^*$.  Thus, $e$ realizes $p$, as required.
 \qed
 \end{proof}

 The following technical Lemma shows that f.g.\ hulls act well under chains.
 
 \begin{Lemma}  \label{fgchain}
 Suppose $(B_n:n\in\omega)$ is an increasing sequence of sets with $B=\bigcup\{B_n:n\in\omega\}$ and suppose $(N_n:n\in\omega)$
 form an elementary chain with $N_0$ an f.g.\ hull over $B_0$
 and each $N_{n+1}$ an f.g.\ hull over $N_n\cup B_{n+1}$.   Assume further that each $N_n$ satisfies $\fg {N_n} {B_n} B$.  Then
 $N^*=\bigcup\{N_n:n\in\omega\}$ is an f.g.\ hull over $B$.
 \end{Lemma}
 
 \begin{proof}  
 First, arguing by induction on $n$, each $N_n$ is an f.g.\ hull over $B_n$ by Lemma~\ref{onestepfg}.  
Next, we show that $N^*$ is f.g.\ full over $B$.   For this, 
 take a finite $\abar$ from $N^*$ and choose a finitely generated type  $p\in S(B\abar)$.
 Suppose $\theta(x,\bbar,\abar)$ generates $p$.  Choose
 $n\in\omega$ so that $\bbar\abar\subseteq N_n$.   Let $p_n\in S(B_n\abar)$ be the restriction of $p$.
 Clearly, $p_n$ is also generated by $\theta(x,\bbar,\abar)$.  Since $N_n$ is f.g.\ full over $B_n$, choose $e\in N_n$ realizing
 $p_n$.   But now, arguing by induction on $\ell\ge n$, the independence and Fact~\ref{useful} imply that $e$ realizes
 the restriction  $p\mr{B_\ell\abar}$.   Thus, $e$ realizes $p$.
 
 Finally, we show that $N^*$ is f.g.\ atomic over $B$.  Choose $\cbar$ from $N^*$.  Say $\cbar\subseteq N_n$.  
 Since $N_n$ is f.g.\ atomic over $B_n$, choose a pp-formula $\delta(\xbar,\bbar)\in \tp^+(\cbar/B_n)$ generating
 $\tp(\cbar/B_n)$.    Arguing as above, using Fact~\ref{useful} and induction,
 we see that $\delta(\xbar,\bbar)$ also generates $\tp(\cbar/B_\ell)$ for every $\ell\ge n$.  Thus, $\delta(\xbar,\bbar)$ generates
 $\tp(\cbar/B)$ as well.
 \qed
 \end{proof}

 \section{Rings with a descending pp-chain and the inverse limit}  \label{fixring}

We now specialize to the theorem's hypothesis.  By Prest's chain-condition criterion for modules, a module is totally transcendental if and only if it satisfies the descending chain condition on pp-definable subgroups; it is enough to consider one-variable pp formulas.  See Prest~\cite[Section 3.1]{Prest} and Ziegler~\cite[Theorem 2.1(2)]{Ziegler}.  Since $R$ is countable, the language of left $R$-modules is countable, and total transcendence is equivalent to $\omega$-stability.

Fix, for the whole of this paper, a countable ring $R$ (with unit) and a sequence  $(\phi_m(x):m\in\omega)$ of pp-definable formulas in one free variable without parameters such that, in the left $R$-module ${}_R R$, 
$$\phi_0(R)> \phi_1(R) > \dots $$ form a strictly descending sequence of subgroups of $R$.
For ease of notation, we sometimes write $G_m$ in place of $\phi_m(R)$.

As each of the subgroups is defined by a pp formula, each  $G_m$ is a right ideal of $R$, but in general it is not a left ideal.

Let $$\F=\varprojlim_{n} R/G_n$$ be the inverse limit.   That is, every $\gamma\in\F$ is a coherent sequence $(\gamma(n):n\in\omega)$, where each $\gamma(n)$ is a coset of $G_n$.  
We endow $\F$ with the inverse limit topology, induced by assuming that
each $R/G_n$ is discrete.  
Since $R$ is countable, the space $\F$ is a Polish space.  For each $n\in\omega$, let $\U_n=\{\gamma\in\F:\gamma(n)=G_n\}$.
As each $G_n$ is an abelian group, $\F$ is an abelian  Polish group with
 $\{\U_n:n\in\omega\}$ a neighborhood base of the identity.  

There is a diagonal map 
$\sigma:R\rightarrow\F$
given by  $r\mapsto (r+G_n)_{n\in\omega}$.
  Note that the image $\sigma[R]$ is dense in $\F$, hence any non-empty, open $\U\subseteq\F$ contains $\sigma(r)$ for at least one $r\in R$.  As the sequence $(G_n:n\in\omega)$ is strictly descending, $\F$ has no isolated points, hence $\F$ has size continuum.

%
%
%
%
%
%
%
%
%

\subsection{The submodules  $L^M$ and  the 2-sided ideal $L^R$}

The descending sequence $G_0>G_1>G_2>\dots$ of pp-definable groups in $\RsubR$ naturally extends to any module $M\preceq\C$.   

\begin{Definition}    \label{L}  {\em  Fix a left $R$-module $M\preceq\C$.
\begin{itemize}
\item  For each $m\in\omega$, let $G_m(M):=\phi_m(M)$, a pp-definable subgroup of $M$.
\item  
Let $G_\infty(M):=\bigcap\{G_m(M):m\in\omega\}$.  

\item  For all $n\in\omega$, let $L_n^M=\{a\in M:$ for every $r\in G_n$, $ra\in G_\infty(M)\}$.
\item  Let $L^M=\bigcup \{L_n^M:n\in\omega\}$.
\end{itemize}
}
\end{Definition}

The verifications of the next two facts follow from the fact that each $G_m(M)$ (and hence $G_\infty(M)$) is a subgroup.

\begin{Fact}  \label{Lfacts}  Let $M$ be any left $R$-module.

\begin{enumerate}
\item  If $k\le n$ then $L_k^M\subseteq L_n^M$.
\item  $L^M$ is a left $R$-module, as is every $L_n^M$.
\end{enumerate}
\end{Fact}

\begin{proof}  (1) is immediate since $k\le n$ implies $G_n\subseteq G_k$.

(2)  As $L^M=\bigcup_{n\in\omega} L_n^M$, it suffices to show that $L_n^M$ is a left $R$-module.  Fix $n\in\omega$ and first choose $a,b\in L_n(M)$.  Then, for each $r\in G_n$, both $ra,rb\in G_\infty(M)$, hence $r(a+b)\in G_\infty(M)$ as well.
So $a+b\in L_n^M$.  Next, suppose $a\in L_n^M$ and $s\in R$.  To show $sa\in L_n^M$, choose any $r\in G_n$.  As $G_n$ is a right ideal $rs\in G_n$ as well.  So 
$$r(sa)=(rs)(a)\in G_{\infty}(M)$$
Thus, $sa\in L_n^M$.
\qed
\end{proof}

\begin{Lemma}   \label{2sided}   $L^R$ is a proper, 2-sided ideal of $R$.
\end{Lemma}

\begin{proof}  By Fact~\ref{Lfacts} we know $L^R$ is a left ideal.  
To see that it is a right ideal, choose $a\in L^R$ and let $s\in R$ be arbitrary.  Say $a\in L_n^R$.
We argue that $as\in L_n^R$ as well.  Indeed, choose $r\in G_n$.  Then $ra\in G_\infty$.  However,
since each $G_m$ is a right $R$-module, $r(as)=(ra)s\in G_\infty$ as well.  

To show properness, we argue that the unit element $1\not\in L^R$.  If it were, then for some $n\in\omega$ we would
have $r\cdot 1=r\in G_\infty$ for every $r\in G_n$.   That is, we would have $G_n=G_\infty$, contradicting our assumption that $G_{n+1}$ is a proper subgroup of $G_n$.
\qed
\end{proof}

In the following Fact, we really only require purity.  

\begin{Fact}  \label{Lpure}   Suppose $M\preceq N\preceq\C$.    Then $L^M=L^N\cap M$.  In fact, $L_n^M=L_n^N\cap M$ for every $n\in\omega$.  
\end{Fact}

\begin{proof}  As the second statement implies the first, fix any $n\in\omega$.  As $G_\infty(M)$ and $G_\infty(N)$ are both defined by the same set of pp formulas, we have $G_\infty(M)=G_\infty(N)\cap M$ since $M$ is pure in $N$.  The result follows easily.
\qed
\end{proof}

%

\begin{Lemma}  \label{freeL}   Let $RX$ be any free $R$-module with an infinite basis set $X$.  Then $L^{RX}=L^RX$, 
i.e., if $a\in RX$ and $a=\sum s_ix_i$ with each $s_i\in R$ and $x_i\in X$, then $a\in L_{RX}$ if and only if, for every $x_i\in X$,
$s_i\in L^R$.  In fact, the analogous result, $L_n^{RX}=L_n^RX$ holds for every $n\in\omega$.
\end{Lemma}  

\begin{proof}  For any $a=\sum s_ix_i\in RX$, then for any $m\in\omega$, since the formula $\phi_m(x)$ defining $G_m(RX)$ is
pp, then Fact~\ref{product} implies that for every $r\in R$, 
$$RX\models\phi_m(ra) \quad \Longleftrightarrow \quad \hbox{$\RsubR\models \phi_m(rs_i)$ for every $x_i\in X$}$$
That $L^{RX}_n=L^R_nX$  for any $n\in\omega$ follows immediately from this, hence $L^{RX}=L^RX$ as well.
\qed
\end{proof}

\begin{Lemma}  \label{directL}  
Suppose $M\preceq\C$ and $M=A\oplus B\oplus D$.
  Then for any $p\in A$ and $q\in B$, if $p+q\in L^{\C}$, then $p\in L^A$ and $q\in L^B$.
\end{Lemma}

\begin{proof}  Let $\ell=p+q$.  Note  that $\ell\in L^M$ by Fact~\ref{Lpure}.  
Let $\pi_{A}:M\rightarrow A$ and $\pi_{B}:M\rightarrow B$ be the canonical projections onto $A,B$, respectively.  
As each is a surjective homomorphism, for every pp formula $\psi(x)$, if $M\models \psi(\ell)$, then Fact~\ref{surjective} implies
$A\models\psi(\pi_{A}(\ell))$ and $B\models\psi(\pi_{B}(\ell))$.   As $\pi_{A(\ell)}=p$ and $\pi_{B}(\ell)=q$, the result follows.
\qed
\end{proof}

\subsection{The coding types $p_\gamma(x,a)$}As notation, for each $\gamma\in\F$, choose representatives $\gamma\<n\>\in R$ so that $\gamma(n)=\gamma\<n\>+G_n$
for each $n\in\omega$.  The choice of representatives will not matter in what follows below.  

\begin{Definition}    {\em  For every $a\in\C$ and $\gamma\in \F$, let 
$$p_\gamma(x,a)=\{\phi_m(x-\gamma\<m\>a):m\in\omega\}$$
}
\end{Definition}

Visibly, each $p_\gamma(x,a)$ is a set of pp formulas.  With Lemma~\ref{consistent} below, $p_\gamma(x,a)$ is satisfiable,
hence is a partial pp type.
The following additivity  facts are easily established.

\begin{Fact}  \label{additivity}  Let $a,a',b,b'\in\C$ and $\gamma,\delta\in\F$ be arbitrary.  Then:

\begin{enumerate}
\item  $
        p_\gamma(b,a)\ \&\ p_\delta(b',a)
        \quad\Rightarrow\quad
        p_{\gamma\pm\delta}(b\pm b',a)$.
\item  $     p_\gamma(b,a)\ \&\ p_\gamma(b',a')
        \quad\Rightarrow\quad
        p_\gamma(b\pm b',a\pm a')$.
\item  For all $r\in R$, 
        $p_{\sigma(r)}(ra,a)$.

\end{enumerate}
\end{Fact}

In what follows, for a fixed $\gamma\in\F$,
we consider larger pp types formed by taking disjoint unions of the types $p_\gamma(x,a)$.

\begin{Definition}  {\em  Let $\abar=(a_0,\dots,a_{n-1})$ be any finite tuple from $\C$.
For any $\gamma\in \F$, $p_\gamma^{\abar}(\xbar,\abar)=\bigcup\{p_\gamma(x_i,a_i):i<n\}$, a set of pp formulas,
 and for any model $M\supseteq\abar$,  let $q^{\abar}_\gamma(\xbar)\in S(M)$ be the complete type generated by $p_\gamma^{\abar}(\xbar,\abar)$.
On occasion, for an infinite tuple $\abar^*$, we write 
$p_\gamma^{\abar^*}(\xbar^*,\abar^*)=\bigcup\{p_\gamma(x_a,a):a\in\abar^*\}$, which is a set of pp formulas in infinitely many variables.
}
\end{Definition}

\begin{Lemma}    \label{consistent}  Fix any $\gamma\in\F$.
\begin{enumerate}
\item For all $\gamma\in \F$ and $a\in\C$, $p_\gamma(x,a)$ is satisfiable.  
\item  For any tuple $\abar$ (finite or infinite) $p^{\abar}_\gamma(\xbar,\abar)$ is satisfiable.  
\item  Suppose $\abar$ is an arbitrary tuple and $\abar'$ is a finite subsequence of $\abar$.
Let $\xbar'$ be the associated finite subsequence of $\xbar$.  If $p^{\abar}_\gamma(\xbar,\abar)$ entails
some pp formula $\delta(\xbar',\cbar)$, then $p^{\abar'}_\gamma(\xbar',\abar')$ also entails $\delta(\xbar',\cbar)$.
\end{enumerate}
\end{Lemma}

\begin{proof}  (1) Fix any $m^*\in\omega$.   
Clearly, the element $1 \in R$ realizes $\phi_m(\gamma \langle m^* \rangle x - \gamma\langle m \rangle x)$ for every $m \leq m^*$. So, by Fact~\ref{1suffices}, $\C \models \forall y(\phi_m(\gamma\langle m^* \rangle y- \gamma \langle m \rangle y))$ for every $m \leq m^*$.
That is, 
 every finite subset of $p_\gamma(x,a)$ is satisfied by $\gamma\<m^*\>a$ for some large enough $m^*$.  
Thus, $p_\gamma(x,a)$  is satisfiable by compactness.

(2) As the types $p_{\gamma}(x_i,a_i)$ and $p_\gamma(x_j,a_j)$ all have distinct free variables, (2) follows immediately from (1).

(3)  Assume that  $p^{\abar}_\gamma(\xbar,\abar)$ entails
some pp formula $\delta(\xbar',\cbar)$ and choose any finite subsequence $\abar'$ from $\abar$.  To see that
$p^{\abar'}_\gamma(\xbar',\abar')$ entails $\delta(\xbar',\cbar)$, choose any $\hbar'$ realizing $p^{\abar'}_\gamma(\xbar',\abar')$.
For each $a\in\abar\setminus \abar'$, by (1) choose some realization $h_a$ of $p_\gamma(x_a,a)$ and let
$\hbar$ be the sequence formed by concatenating $\hbar'$ with each of the $h_a$'s (in the correct order).  Then
$\hbar$ realizes $p_\gamma^{\abar}(\xbar,\abar)$, hence $\delta(\hbar',\cbar)$ by assumption.
\qed
\end{proof}

Each of the following statements is verified by checking  that each $H_i$ is closed under subtraction.  

\begin{Lemma}  \label{subgroups}   Suppose $\abar,e$ are fixed, $\phi(y,\xbar,\zbar)$ is a pp-formula with $\lg(\xbar)=\lg(\abar)$, $\lg(\zbar)=\lg(\dbar)$, 
and $\delta\in F$.   Then each of the following sets is either empty or is a subgroup.

\begin{enumerate}
\item  $H_1=\{(\gamma,\dbar)\in \F\times\C^{\lg(\zbar)}:p^{\abar}_\gamma(\xbar,\abar)\models \exists y \phi(y,\xbar,\dbar)\}$
\item  $H_2:=H_1\cap \{(\gamma,\dbar)\in \F\times\C^{\lg(\zbar)}:p^{\abar}_\gamma(\xbar,\abar)\wedge \phi(y,\xbar,\dbar)\models p_\gamma(y,e)\}$
\item $H_3:=H_1\cap \{(\gamma,\dbar)\in \F\times\C^{\lg(\zbar)}:p^{\abar}_\gamma(\xbar,\abar)\wedge \phi(y,\xbar,\dbar)\models p_\delta(0,y)\}$
\item $H_4:=(H_1\times\C)\cap \{(\gamma,\dbar c)\in \F\times\C^{\lg(\zbar)+1}:p^{\abar}_\gamma(\xbar,\abar)\wedge \phi(y,\xbar,\dbar)\models p_\delta(y,c)\}$.
\end{enumerate}
\end{Lemma}


The following is a basic result about groups.

\begin{Fact}  Suppose $C,D\subseteq G$ are cosets of the subgroups $H,K$, respectively.
Then $C\subseteq D$ if and only if $C\cap D\neq\emptyset$ and $H\le K$.
\end{Fact}

This translates immediately in the context of pp-formulas.

\begin{Fact}  \label{ppcosets}   Let $\phi(\xbar,\ybar)$ and $\psi(\xbar,\zbar)$ be pp-formulas without parameters.
\begin{enumerate}
\item  For every $\abar$, $\phi(\C,\abar)$ is either empty or is a coset of $\phi(\C,\0bar)$.
\item  If, for some $\abar,\bbar$,  $\phi(\C,\abar)\neq\emptyset$ and 
$\phi(\xbar,\abar)$ entails $\psi(\xbar,\bbar)$, then the subgroups $\phi(\C,\0bar)\le \psi(\C,\0bar)$.
\item  Conversely, if $\phi(\C,\0bar)\le \psi(\C,\0bar)$, then for any $\abar,\bbar$, if $\exists\xbar(\phi(\xbar,\abar)\wedge\psi(\xbar,\bbar))$, then $\phi(\xbar,\abar)$ entails $\psi(\xbar,\bbar)$.
\end{enumerate}
\end{Fact}

\begin{Lemma}  \label{typesmodL}  Let $M$ be any left $R$-module, with $a\in M$ and $\gamma\in\F$.  If
$p_\gamma(x,a)$ is realized in $M$, then $p_\gamma(x,b)$ is also realized in $M$ for every $b\in a+L^M$.
In particular, if $e\in M$ realizes $p_\gamma(x,a)$, then for every $\ell\in L^M$ there is some $s\in R$ such that
$p_\gamma(e+s\ell,a+\ell)$ holds.
\end{Lemma}

\begin{proof}  The second statement implies the first, so assume $e\in M$, $\ell\in L^M$, and $p_\gamma(e,a)$ holds. 
Suppose $\ell\in L_n^M$.   As $\sigma[R]$ is dense in $\F$, choose $s\in R$ so that $\gamma-\sigma(s)\in \U_n$.  

\medskip
\noindent{\bf Claim.}  For all $m\in\omega$, $(\gamma\<m\>-s)\ell\in G_m(M)$.

\begin{proof}  There are two cases.  First, assume $m\le n$.  As $\gamma\<m\>-\gamma\<n\>\in G_m$, $\gamma\<n\>-s\in G_n$, and $G_n\le G_m$, we have
 $(\gamma\<m\>-s)\in G_m$.   As the map $r\mapsto r\ell$ is a homomorphism between $R$ and $M$, we have $(\gamma\<m\>-s)\ell\in G_m(M)$.
On the other hand, if $m\ge n$, then at least $(\gamma\<m\>-s)\in G_n$.  But then, as $\ell\in L_n(M)$ we have
$(\gamma\<m\>-s)\ell\in G_\infty(M)\subseteq G_m(M)$.
\end{proof}


Restating the  Claim,  $\gamma\<m\>\ell\equiv s\ell$ (mod $G_m(M)$) for every $m\in \omega$.  But also,
$p_\gamma(e,a)$ holding implies $\gamma\<m\>a\equiv e$ (mod $G_m(M)$) for each $m\in \omega$.
Combining these gives $\gamma\<m\>(a+\ell)\equiv e+s\ell$ (mod $G_m(M)$) for every $m$, i.e., 
$p_\gamma(e+s\ell,a+\ell)$ holds.
\qed
\end{proof}

 \section{From tagged $R$-modules to $R$-modules}   \label{transfer}
 
 In this section, we start with a countable $R$-module $M$ with a countable set of distinguished 
 submodules $\Abar=(A_j:j\in J)$ such that each $A_j$ is a direct summand in $M$.   Our goal will be to define a `suitably generic' doubly indexed sequence $\gammabar=(\gamma_{j,m}:j\in J,m\in\omega)$ from $\F$ and build a  countable $R$-module
 $M(\Abar,\gammabar)\succeq M$ that `encodes the submodules $\Abar$' in terms of realizing and omitting pp types.  This is made precise in Theorem~\ref{onegamma}.    
 
 As notation, assume that $\abar_j^*$ is an enumeration of $A_j$ and let $\xbar^*$ be a sequence of distinct variable symbols
 of the same order type.   
 
 We set some notation.  Suppose a sequence $\gammabar=(\gamma_{j,m}:j\in J,m\in\omega)$  has been chosen.

\begin{itemize}
\item  $\xtilde=(x_{j,m,a}:j\in J, m\in\omega, a\in \abar_j^*)$ is a triply-indexed sequence of distinct variable symbols.  
\item  Let $p^{\Abar}(\xtilde)=\bigcup\{p_{\gamma_{j,m}}(x_{j,m,a},a):j\in J,m\in\omega,a\in \abar_j^*\}$.
\item Let $q^{\Abar}(\xtilde)\in S(M)$ be the complete type generated by $p^{\Abar}(\xtilde)$.
\item  Let $H$ denote a realization of $q^{\Abar}$ and 
\item  Let $M(\Abar,\gammabar)$ be an f.g.\ hull of $MH$.
%
\end{itemize}

The following Theorem clarifies what we mean by `$\gammabar$ is sufficiently generic'. 


\begin{Theorem}  \label{onegamma}  Let $M\preceq\C$ be a countable $R$-module with a countable set of distinguished 
 submodules $\Abar=(A_j:j\in J)$ such that each $A_j$ is a direct summand in $M$.
 Then there is a sequence $\gammabar=(\gamma_{j,m}:j\in J,m\in\omega)$ from $\F$ such that,
 for every $j\in J$ and $b\in M(\Abar,\gammabar)$,
 $$\hbox{$M(\Abar,\gammabar)$ realizes $p_{\gamma_{j,m}}(x,b)$ for every $m\in\omega \Longleftrightarrow
 b\in A_j+L^{M(\Abar,\gammabar)}$}$$
 \end{Theorem}

 Before launching into a proof of Theorem~\ref{onegamma}, we introduce a new notion.
 
 \begin{Definition}  {\em   Suppose $\psi(w,\zbar)$ is a 0-definable pp formula and $b,\abar$ are from $\C$.
 We say {\em  $b$ is linked to $\abar$ via $(n,\psi(w,\zbar))$} if, for all $r\in G_n$
 \begin{enumerate}
 \item  $p_{\sigma(r)}^{\abar}(\zbar,\abar)\models \exists w \psi(w,\zbar)$; and
 \item  $p_{\sigma(r)}^{\abar}(\zbar,\abar)  \wedge \psi(w,\zbar) \models p_{\sigma(r)}(w,b)$.
 \end{enumerate}
 
 For $A$ any set, we say {\em $b$ is linked to $A$} if there is some finite $\abar$ from $A$ and some $(n,\psi(w,\zbar))$
 with $b$ linked to $\abar$ via $(n,\psi(w,\zbar))$.
 }
 \end{Definition}

 Trivially, if $b\in A$, then $b$ is linked to $A$, taking $\psi(w,z)$ to be `$w=z$'.  
 As is the case with realizing $p_\gamma$'s, this notion is closed under adding elements of $L^{\C}$.
 
 \begin{Lemma}  If $b$ is linked to $A$, then so is $b+\ell$ for any $\ell\in L^{\C}$.
 \end{Lemma}
 
 \begin{proof}  Suppose $\abar\subseteq A$ and $b$ is linked to $\abar$ via $(n,\psi(w,\zbar))$.   Also, suppose $\ell\in L_m^{\C}$.
 Let $k=\max\{n,m\}$ and we show $b+\ell$ is linked to $\abar$ via $(k,\psi(w,\zbar))$.  To see this, choose any $r\in G_k$.  Choose any $\hbar$ realizing 
 $p^{\abar}_{\sigma(r)}(x,\abar)$ and choose any $d$ such that $\psi(d,\hbar)$.  Then
 $p_{\sigma(r)}(d,b)$ holds.   But also, $\ell\in L_k^{\C}$ implies $p_{\sigma(r)}(0,\ell)$ holds.  Thus, $p_{\sigma(r)}(d,b+\ell)$ holds, as required.
 \qed
 \end{proof}
 
 Thus, for any $M$ with distinguished submodules $\Abar=(A_j:j\in J)$, if $b\in A_j+L^M$, then $b$ is linked to $A_j$.
 The following definition describes the converse.

\begin{Definition} {\em  Suppose $M\preceq\C$ and $\Abar=(A_j:j\in J)$ are distinguished submodules of $M$.  We say that
{\em  $M$ is correct for $\Abar$} if, for all $j\in J$ and for all $b\in M$, if $b$ is linked to $A_j$, then $b\in A_j+L^M$.
}
\end{Definition}

\begin{Lemma}  \label{correct0}  
Assume $M\preceq \C$ and every $A_j$ is a direct summand of $M$.  Then $M$ is correct for $\Abar$.
\end{Lemma}

\begin{proof}  Fix $j\in J$ and choose a submodule $K$ such that $M=A_j\oplus K$.  Let $\pi_1,\pi_2$ denote the projection maps onto $A_j,K$, respectively.  Choose any $b\in M$ that is linked to $A_j$.    As $b=\pi_1(b)+\pi_2(b)$ with $\pi_1(b)\in A_j$,
it suffices to show that $\pi_2(b)\in L^M$.    For this, suppose $b$ is linked to $\abar\subseteq A_j$ via $(n,\psi(w,\zbar))$.
We argue that $\pi_2(b)\in L_n^M$.  To see this, choose $r\in G_n$.  By Fact~\ref{additivity}(3), the tuple $r\abar$ realizes
$p_{\sigma(r)}^{\abar}(\xbar,\abar)$, and  since $\abar\subseteq A_j$, $r\abar\subseteq A_j$ as well.  Thus, there is some $d\in\C$ such that
$\psi(d,r\abar)$ holds.   As $M\preceq \C$ and 
since $A_j$ is pure in $M$, there is some $d'\in A_j$ such that $\psi(d',r\abar)$ holds.  By the second entailment of being linked,
we have $p_{\sigma(r)}(d',b)$ holding.   Now, applying $\pi_2$, which is a surjective homomorphism onto $K$, we have
$p_{\sigma(r)}(\pi_2(d'),\pi_2(b))$ holding.   Since $d'\in A_j$,  this gives $p_{\sigma(r)}(0,\pi_2(b))$.   As this holds for every $r\in G_n$,
$\pi_2(b)\in L^M$.
\qed
\end{proof}

\begin{Notation} {\em  Suppose $M\preceq\C$ is countable,  $A\subseteq M$ is a submodule, and $\gamma\in \F$.  Fix an enumeration $\abar^*$ of $A$ and a matching sequence $\xbar^*$ of distinct variable symbols.  
\begin{itemize}
\item  $p^A_\gamma(\xbar^*,\abar^*):=\bigcup\{p_\gamma(x_a,a):a\in A\}$.
\item  $q^A_\gamma(\xbar^*)$ denotes the complete type in $S_{\xbar^*}(M)$ generated by $p^A_\gamma(\xbar^*,\abar^*)$.
\item  $\hbar^*_\gamma$ denotes an arbitrary realization of $q^A_\gamma$.
\item    $M(A,\gamma)$ denotes a countable f.g.\ hull over $M\hbar^*_\gamma$.
\end{itemize}
}
\end{Notation}

Note that since $q^A_\gamma(\xbar^*)$ is a complete type and f.g.\ hull over $M\hbar^*_\gamma$
is unique up to isomorphism over $M\hbar^*_\gamma$, the isomorphism type of $M(A,\gamma)$ is independent of the choices made.  
The following Proposition describes the `generic behavior.'  Its proof is the content of Subsection~\ref{4.1}.

\begin{Proposition} \label{genericfacts}
 Suppose $M\preceq\C$ is countable and $\Abar=(A_j:j\in J)$ are countably many distinguished submodules.
\begin{enumerate}[label=\Roman*.]
\item  For every $\delta\in\F$, $b\in M$, $j\in J$, if $p_\delta(x,b)$ is omitted in $M$, then
$$\{\gamma\in\F:\hbox{$M(A_j,\gamma)$ omits $p_\delta(x,b)$}\} \ \hbox{is comeager}$$
\item  If $j\in J$ then
$$\{\gamma\in\F:\hbox{$M(A_j,\gamma)$ omits $p_\gamma(x,b)$ for every $b\in M$ not linked to $A_j$}\} \ \hbox{is comeager}$$
\item If $M$ is correct for $\Abar$, then for every $j\in J$
$$\{\gamma\in\F:\hbox{$M(A_j,\gamma)$ is correct for $\Abar$}\} \ \hbox{is comeager}$$
\end{enumerate}
\end{Proposition}

\noindent{\bf Proof of Theorem~\ref{onegamma} given Proposition~\ref{genericfacts}.}

Our proof will be rather indirect.   We will simultaneously construct an elementary $\omega$-chain $M_k$ of models and
construct our sequence $\gammabar$ in $\omega$ steps.  We will then argue that $N^*$, the union of the chain
is isomorphic to $M(\Abar,\gammabar)$ as defined above.

Fix a bijection $g:\omega\rightarrow J\times\omega$.
Let $M_0:=M$.   By Lemma~\ref{correct0}, $M$ is correct for $\Abar$.
Now suppose $M_k$ has been defined and is correct for $\Abar$, and additionally assume $\gamma_{g(\ell)}$ has been chosen for all $\ell<k$.  Suppose $g(k)=(j,m)$.  In light of Proposition~\ref{genericfacts}, we can choose $\gamma^*$ such that:

\begin{enumerate}
\item  For every $b\in M_k$ and for every $\ell< k$, if $p_{\gamma_{g(\ell)}}(x,b)$ is omitted in $M_k$, then $p_{\gamma_{g(\ell)}}(x,b)$ is omitted
in $M_k(A_j,\gamma^*)$;
\item  $M_k(A_j,\gamma^*)$ is correct for $\Abar$;
\item  $M_k(A_j,\gamma^*)$ omits $p_{\gamma^*}(x,b)$ for every $b\in M_k$ that is not linked to $A_j$.
\end{enumerate}
Choose a realization $\hbar_k^*$ of $q^{A_j}_{\gamma^*}(\xbar^*)$ and let
$M_{k+1}:=M_k(A_j,\gamma^*)$ be a countable, f.g.\ hull over $M_k\hbar_k^*$.
Also, put $\gamma_{g(k)}=\gamma_{j,m}:=\gamma^*$.  

Finally, put $N^*:=\bigcup\{M_k:k\in\omega\}$ and $\gammabar=(\gamma_{g(k)}:k\in\omega)=(\gamma_{j,m}:j\in J,m\in\omega)$.
As every $M_k$ is correct for $\Abar$, it is easily checked that $N^*$ is as well.  
We check that $N^*$ satisfies the omitting types condition of the Theorem.  
    Fix $j\in J$.   First,  if $a\in A_j$, then for every $m\in\omega$, the element $h_a\in\hbar^*_{j,m}$ realizes 
    $p_{\gamma_{j,m}}(x,a)$.
Thus, by Lemma~\ref{typesmodL}, if $b\in A_j+L^{N^*}$, then $p_{\gamma_{j,m}}(x,b)$ is realized in $N^*$.

Conversely,  suppose $b\in N^*\setminus(A_j+L^{N^*})$.  Since $N^*$ is correct for $\Abar$, this implies
$b$ is not linked to $A_j$.  
Choose $k$ sufficiently large  so that $b\in M_k$ and $g(k)=(j,m)$ for some $m\in\omega$.  
Then by (3), $p_{\gamma_{g(k)}}(x,b)$ is omitted in $M_{k+1}$.   Because of (1), it follows that $p_{\gamma_{g(k)}}(x,b)$ 
remains omitted in $M_\ell$ for every $\ell\ge k+1$, hence $p_{\gamma_{j,m}}(x,b)$ is omitted in $N^*$.


To finish the proof of the Theorem, it suffices to prove that $N^*$ is an f.g.\ hull over $MH$, where $H=\bigcup\{\hbar^*_k:k\in\omega\}$.  For this we use Lemma~\ref{fgchain}.  For each $k$, let
$B_k:=M_0\cup\bigcup\{\hbar^*_i:i\le k\}$ and let $N_k:=M_{k+1}$.   Then
\begin{itemize}
\item  $N_0=M_1$ is an f.g.\ hull over $B_0=M_0\hbar^*_0$;
\item  $N_{k+1}=M_{k+2}$ is an f.g.\ hull over $M_{k+1}\hbar^*_{k+1}=N_kB_{k+1}$; and
\item  $B=\bigcup_k B_k=MH$ and $\bigcup_k N_k=N^*$.
\end{itemize}

Additionally,  each $\hbar^*_k$ realizes a type over $M_k$ that is generated by
a pp-partial type over $M_0$.  Thus, it follows by non-forking symmetry that $\fg {M_k} {M_0} {\hbar^*_k}$.   As this holds for every $k$, using symmetry and transitivity it follows by induction 
that $\fg {N_k} {B_k} B$ for every $k\in\omega$.     Finally,  By induction on finite subtuples and Proposition~\ref{main}, the type of $H$ over $M$ is generated by $p^{\Abar}$,
so $H$ realizes the complete type $q^{\Abar}$. 
Thus,  $N^*$ is an f.g.\ hull over $MH$ by Lemma~\ref{fgchain}.
\qed

\subsection{Proof of Proposition~\ref{genericfacts}}  \label{4.1}

Here we prove the three density conditions.   Each of these will be proved by an application of Pettis' theorem, asserting that any non-meager coset of a Borel subgroup of a Polish group is open.    The Borelness of the cosets $\Gamma_\partial$ and $\Gamma_{\partial'}$ follows from this fact,  whose verification is  left to the reader.

\begin{Fact}  \label{Borelness}  Fix an enumeration of the ring $R$ and, for each $m\in\omega$ choose $\gamma\<m\>\in R$ to be least in its $R/G_m$-class.  Then, for fixed parameters, fixed pp formulas, and fixed $\delta\in\F$,
\begin{itemize}
\item  $\{\gamma\in\F: p_\gamma^{\abar}(\xbar,\abar)\models \theta(\xbar,\cbar)\}$ is open;
\item  The sets $\{\gamma\in\F: p_\gamma^{\abar}(\xbar,\abar)\wedge\phi(y,\xbar,\dbar)\models p_\delta(y,c)\}$ and
$\{\gamma\in\F: p_\gamma^{\abar}(\xbar,\abar)\wedge\phi(y,\xbar,\dbar)\models p_\gamma(y,c)\}$ are $G_\delta$;
\item  The extended template $\Gamma_{\partial'}$ in {\bf (III)} below is Borel since the intersection over $r\in G_n$ is countable.
\end{itemize}
\end{Fact}  

We begin with a lemma describing how elements of
$M(A,\gamma)$ can be described by sets of pp formulas over $M$.

\begin{Lemma}  \label{p-part}  Suppose $M$ is countable, $A\subseteq M$, and $\gamma\in \F$.  Suppose $M(A,\gamma)$ is the f.g.\ hull
over $M\hbar$, where $\hbar$ realizes $q^A_\gamma(\xbar)\in S(M)$.   Then for any $b\in M(A,\gamma)$,
there is a finite tuple $\abar$ from $A$, a pp formula $\phi(y,\xbar,\zbar)$, and a tuple $\dbar$ from $M$ so that
$p^{\abar}_\gamma(\xbar,\abar)\wedge\phi(y,\xbar,\dbar)$ entails the pp type $\tp^+(b/M)$.
\end{Lemma}

\begin{proof}  As $\tp(b/M\hbar)$ is finitely generated, choose $\phi(y,\xbar,\dbar)$ with $\dbar$ from $M$ so that
$\phi(y,\hbar,\dbar)$ generates $\tp(b/M\hbar)$.   
Only finitely many coordinates of $\hbar$ occur in this generator.
Replace $\hbar$ by that finite subtuple, and let $\abar$ be the
corresponding finite tuple from $A$.  By Lemma~4.10(3), the
restriction \(\tp(\hbar/M)\) is generated by
$p_\gamma^{\abar}(x,\abar)$.

Choose any $\delta(y,\mbar)\in \tp^+(b/M)$.   As $\phi(y,\hbar,\dbar)$
entails $\delta(y,\mbar)$, by Fact~\ref{ppcosets}(2) we have $\phi(\C,\0bar,\0bar)\le \delta(\C,\0bar)$.   
Also, the pp-formula $\exists y(\phi(y,\hbar,\dbar)\wedge\delta(y,\mbar))$ is witnessed by $b$.
Thus, $$\exists y(\phi(y,\xbar,\dbar)\wedge\delta(y,\mbar))\in\tp^+(\hbar/M).$$
Since 
$p^{\abar}_\gamma(\xbar,\abar)$ generates $\tp(\hbar/M)$, it entails
$\exists y(\phi(y,\xbar,\dbar)\wedge\delta(y,\mbar))$.
Thus, $p^{\abar}_\gamma(\xbar,\abar)\wedge\phi(y,\xbar,\dbar)$ entails $\delta(y,\mbar)$ by Fact~\ref{ppcosets}(3).
As this entailment holds for every $\delta(y,\mbar)\in\tp^+(b/M)$, this suffices.
\end{proof}

\noindent {\bf (I)}   Fix $A_j\in\Abar$ and some $e\in M$ so that $p_\delta(y,e)$ is omitted in $M$.
  As $M$ is countable, it suffices to prove that
$$Y=\{\gamma\in\F:M(A_j,\gamma)\ \hbox{realizes $p_\delta(y,e)$}\}$$
is meager.  
A {\em template} $\partial=(\abar,\dbar,\phi(y,\xbar,\zbar))$ consists of finite sets $\abar\subseteq A_j$, $\dbar\subseteq M$ and a pp-formula $\phi(y,\xbar,\zbar)$  without parameters.  
Associated to each template $\partial$, put
$$\Gamma_\partial:=\{\gamma\in\F: p_\gamma^{\abar}(\xbar,\abar)\models \exists y\phi(y,\xbar,\dbar)\ \hbox{and} \
p_\gamma^{\abar}(\xbar,\abar)\wedge\phi(y,\xbar,\dbar)\models p_\delta(y,e)\}$$
Note that there are only countably many templates and by Lemma~\ref{p-part}, $Y\subseteq\bigcup\{\Gamma_\partial:\partial\}$.
Thus, to prove that $Y$ is meager, it suffices to prove the following Claim.

\medskip
\noindent{\bf Claim.}  For all templates $\partial$, $\Gamma_\partial$ is meager.

\begin{proof}  By way of contradiction, choose a template $\partial$ for which $\Gamma_\partial$ is non-meager.
Note that $$\Gamma_\partial=((H_1\times\C)\cap H_4)_{\dbar e}$$ in the notation of Lemma~\ref{subgroups}.   As we are assuming it is non-meager, hence non-empty, $\Gamma_\partial$ is a 
coset of a subgroup $K\le\F$.   It follows from Pettis' theorem that $\Gamma_\partial$ is open.   
As $R$ is dense in $\F$, choose $r\in R$ such that $\sigma(r)\in \Gamma_\partial$.  Since $\abar\subseteq A_j\subseteq M$,
$r\abar=(ra_1,\dots,ra_n)\subseteq M$ and $p_{\sigma(r)}(r\abar,\abar)$ holds by Fact~\ref{additivity}(3).  As $\sigma(r)\in \Gamma_\partial$,
$\exists y\phi(y,r\abar,\dbar)$, so as $M\preceq\C$, choose $c\in M$ such that $\phi(c,r\abar,\dbar)$ holds.  
But then, by the second clause of $\Gamma_\partial$, $p_\delta(c,e)$ holds, contradicting our assumption that $p_\delta(y,e)$ is omitted in $M$.
\qed
\end{proof}

\noindent{\bf (II)}  Fix $A_j\in\Abar$ and some $e\in M$ such that $e$ is not linked to $A_j$.
 Let $Y=\{\gamma\in \F:M(A_j,\gamma)$ realizes $p_\gamma(y,e)\}$.    
As above, choose a template $\partial=(\abar,\dbar,\phi(y,\xbar,\zbar))$ and here let
$$\Gamma_\partial:=\{\gamma\in\F:p_\gamma^{\abar}(\xbar,\abar)\models \exists y\phi(y,\xbar,\dbar)\ \hbox{and} \
p_\gamma^{\abar}(\xbar,\abar)\wedge\phi(y,\xbar,\dbar)\models p_\gamma(y,e)\}$$
As above, it suffices to prove the following claim.

\medskip
\noindent{\bf Claim.}  For all templates $\partial$, $\Gamma_\partial$ is meager.

\begin{proof}  By way of contradiction, assume that $\Gamma_\partial$ is non-meager.  Then, in the notation of
Lemma~\ref{subgroups}, $\Gamma_\partial=(H_1\cap H_2)_{\dbar}$.   As it is non-empty, it is a coset of 
$$K_\partial=\{\gamma\in\F:p_\gamma^{\abar}(\xbar,\abar)\models \exists y\phi(y,\xbar,\0bar)\ \hbox{and} \
p_\gamma^{\abar}(\xbar,\abar)\wedge\phi(y,\xbar,\0bar)\models p_\gamma(y,e)\}$$
By Pettis' theorem, $K_\partial$ is open.   Choose $n$ for which $\U_n\subseteq K_\partial$
and let $\psi(y,\xbar):=\phi(y,\xbar,\0bar)$.  Then for every $r\in G_n$, since $\sigma(r)\in K_\partial$, we get
that $e$ is linked to $\abar$ via $(n,\psi(y,\xbar))$, contradicting our assumption on $e$.
\qed
\end{proof}

\noindent{\bf (III)}  Assume that $M$ is correct and $j\in J$.   Call a $\gamma\in\F$ {\em bad} if there is some $b\in M(A_j,\gamma)$ and
some $A_\ell\in \Abar$ such that $b$ is linked to $A_\ell$, but $b\not\in A_\ell+L^{\C}$.  
To subdivide the bad elements, define an extended template 
$$\partial'=(\abar,\dbar,\phi(y,\xbar,\zbar),\ell,\abar',n,\psi(w,\zbar))$$
where $\abar,\dbar,\phi(y,\xbar,\zbar)$ are as above, $\ell\in J$, $\abar'\subseteq A_\ell$ is finite, $n\in\omega$, and 
$\psi(w,\zbar)$ is a pp-formula 
such that $p_{\sigma(r)}^{\abar'}(\zbar,\abar')\models \exists w\psi(w,\zbar)$ for every $r\in G_n$.  

For each extended template $\partial'$, let $\Gamma_{\partial'}$ be the set of all $\gamma\in \F$ satisfying

\begin{equation}  \label{1}
p^{\abar}_\gamma(\xbar,\abar)\models \exists y \phi(y,\xbar,\dbar)
\end{equation}
and, for every $r\in G_n$,
\begin{equation}  \label{2}
p^{\abar}_\gamma(\xbar,\abar)\wedge\phi(y,\xbar,\dbar)\wedge p_{\sigma(r)}^{\abar'}(\zbar,\abar')\wedge \psi(w,\zbar)\models
p_{\sigma(r)}(w,y)
\end{equation}  

\medskip
\noindent{\bf Claim 1.}  If $\gamma\in \F$ is bad, then $\gamma\in \Gamma_{\partial'}$ for some extended template $\partial'$.

\begin{proof}  Suppose $b\in M(A_j,\gamma)$ witnesses that $\gamma$ is bad.  Choose a template 
$\partial=(\abar,\dbar,\phi(y,\xbar,\zbar))$ as above and suppose $b$ is linked to $\abar'\subseteq A_\ell$ via $(n,\psi(w,\zbar))$.
We show $\gamma\in \Gamma_{\partial'}$ for  this extended template.  Clause (1) is as above.  For (2), fix $r\in G_n$ and choose any $m\in\omega$.   We use  Fact~\ref{ppcosets}(3) to show that the antecedent of (2) implies $\phi_m(w-ry)$.   As $m$ is arbitrary, this suffices.

Since $b$ is linked to $\abar'$ via $(n,\psi(w,\zbar))$ and compactness, there
is $\theta_{r,m}(\zbar,w)$ a finite conjunction from $p^{\abar'}_{\sigma(r)}(\zbar,\abar')\wedge\psi(w,\zbar)$ so that
$\theta_{r,m}(\zbar,w)\models \phi_m(w-rb)$.  Let
$$\eta_{r,m}(y):=\exists\zbar w (\theta_{r,m}(\zbar,w)\wedge\phi_m(w-ry))$$
Since $\eta_{r,m}(y)\in\tp^+(b/M)$, it is entailed by $p^{\abar}_\gamma(\xbar,\abar)\wedge\phi(y,\xbar,\dbar)$.
Thus, the antecedent of Clause (2) entails $\phi_m(w-ry)$ by Fact~\ref{ppcosets}(3).
\qed
\end{proof}


\medskip
\noindent{\bf Claim 2.}  For all extended templates $\partial'$, if $\Gamma_{\partial'}$ is non-meager,
then  for every $\gamma\in \Gamma_{\partial'}$ and every $\ubar_\gamma,t_\gamma$ realizing 
$p^{\abar}_\gamma(\xbar,\abar)\wedge\phi(y,\xbar,\dbar)$, we have $t_\gamma\in A_\ell+L^{\C}$.

Note that, given the  two Claims, every $\gamma\in\F$ that is bad is contained in some $\Gamma_{\partial'}$ that is meager.
As there are only countably many extended templates, this implies that the set of bad $\gamma$'s is meager, proving (3).
Thus, it suffices to prove Claim 2. 

\begin{proof}  Fix $\partial'$ for which $\Gamma_{\partial'}$ is non-meager. For each $r\in G_n$, let
\begin{quotation}\noindent
$K_r:=\{u\in \F: p_u^{\abar}(\xbar,\abar)\models \exists y \phi(y,\xbar,\0bar)$ and
$p_u^{\abar}(\xbar,\abar)\wedge\phi(y,\xbar,\0bar)\models p_{\sigma(r)}(0,y)\}$ and let
$K=\bigcap_{r\in G_n} K_r.$
\end{quotation}

It is easily checked that if $\gamma_1,\gamma_2\in \Gamma_{\partial'}$, then $\gamma_1-\gamma_2\in K$, and dually,
if $\gamma\in \Gamma_{\partial'}$ and $u\in K$, then $\gamma+u\in\Gamma_{\partial'}$.
Since $\Gamma_{\partial'}$ was assumed to be non-meager,  $K$ is non-empty.   By Lemma~\ref{subgroups} each $K_r$ and hence 
$K\le\F$ is a subgroup and $\Gamma_{\partial'}$ is a coset of $K$.  As $K$ is Borel and 
we are assuming $\Gamma_{\partial'}$ is non-meager, $\Gamma_{\partial'}$ is open in $\F$.
Choose any $\gamma\in \Gamma_{\partial'}$ and choose $s\in R$ such that $\gamma-\sigma(s)\in K$.   
As $\Gamma_{\partial'}$ is a $K$-coset, $\sigma(s)\in \Gamma_{\partial'}$
as well.
Since $\abar\subseteq M$, so is the tuple $s\abar$.   As $\exists y\phi(y,s\abar,\dbar)$ holds, choose some $m^*\in M$ 
such that $\phi(m^*,s\abar,\dbar)$.

Choose any $\ubar_\gamma t_\gamma$ realizing $p^{\abar}_\gamma(\xbar,\abar)\wedge \phi(y,\xbar,\dbar)$.   We show $t_\gamma\in A_\ell+L^{\C}$ in two steps.

\medskip
\noindent
{\bf Subclaim a).}  $t_\gamma-m^*\in L_n^{\C}$.

\begin{proof}  For each $r\in G_n$, choose $\hbar_r,c_r$ realizing $p_{\sigma(r)}^{\abar'}(\zbar,\abar')$ and $\psi(c_r,\hbar_r)$.
As both $\gamma,\sigma(s)\in \Gamma_{\partial'}$, the second clause of the definition gives, $p_{\sigma(r)}(c_r,t_\gamma)$
and $p_{\sigma(r)}(c_r,m^*)$ for each $r\in G_n$.  Subtracting and Fact~\ref{additivity} give $p_{\sigma(r)}(0,t_\gamma-m^*)$ for every $r\in G_n$,
which means $t_\gamma-m^*\in L_n^{\C}$ by definition of $L_n$.
\qed
\end{proof}

\medskip
\noindent
{\bf Subclaim b).}  $m^*\in A_\ell+L^M$.

\begin{proof}  Since $s\abar,m^*$ realize $p_{\sigma(s)}^{\abar}(\xbar,\abar)\wedge\phi(y,\xbar,\dbar)$,
$\sigma(s)\in\Gamma_{\partial'}$ implies that $m^*$ is linked to $\abar'$ via $(n,\psi)$.  In particular,
$m^*$ is linked to $A_\ell$.   As $M$ is correct and $m^*\in M$, we conclude that $m^*\in A_\ell+L^M$.
\qed
\end{proof}

Given the two Subclaims, it follows immediately that $t_\gamma\in A_\ell+L^{\C}$, as required.
\qed

\section{Two universal objects}

In this section we define two basic objects,  a countable, tagged $R$-module $F^*$ and a countable $R$-module $N^*$.
Per se, neither of these objects code any information, but our eventual proof of Borel completeness will use suitable substructures of $F^*$ and elementary substructures
of $N^*$.   


Throughout this paper, fix a distinguished, countably infinite set $X$
and let $RX$ denote the (countable) direct sum $\bigoplus_{x\in X}Rx$.  

\subsection{The free-like tagged $R$-module $F^*$}

The construction of $F^*$ mimics a construction in Subsection 3.1 of \cite{LU}.

Let $\{Y_n:n\in\omega\}$ be a set of disjoint, countably infinite sets, each disjoint from $X$, and let
$Z=X\cup\bigcup\{Y_n:n\in\omega\}$.  

For each $n\in\omega$, fix a bijection $\rho_n:Y_n\rightarrow RX$ and let $\rhobar_n:RY_n\rightarrow RX$ be the (unique) $R$-linear map extending $\rho_n$.  Working in the module $RZ$, let 
$W_n:=\{y-\rho_n(y):y\in Y_n\}$ and let $RW_n$ be the submodule freely generated by $W_n$.

\begin{Fact}  \label{Wfacts}  Fix $n\in\omega$.  
\begin{enumerate}
\item  $RW_n=\{b-\rhobar_n(b):b\in RY_n\}$.
\item  $RX\cap RW_n=(0)$, and moreover
\item  As submodules of $RZ$, $RX\oplus RW_n=RX\oplus RY_n$.
\end{enumerate}
\end{Fact}

Going forward, let $J=\{*\}\cup(\omega\times\{0,1\})$ be an index set.

\begin{Definition}  {\em  $F^*$ is the tagged $R$-module $(RZ,A_j:j\in J)$ with universe $RZ$, where
$A_*:=RX$, $A_{(n,0)}=RY_n$, and $A_{(n,1)}=RW_n$.
}
\end{Definition}

Visibly, the universe of $F^*$ is a free $R$-module.  It is easily seen that each $A_{(n,0)}$ is both a free module with a free complement, i.e., $RZ=A_{(n,0)}\oplus C_n$ for some free $R$-module $C_n$ (in this case $C_n=R(Z\setminus Y_n))$.
However, by Fact~\ref{Wfacts}(3), $A_{(n,1)}$ is also freely complemented by $C_n$, i.e.,
$RZ=RW_n\oplus C_n$.

\subsection{The universal $R$-module $N^*$}

Here, we apply Theorem~\ref{onegamma} to $F^*=(RZ,\Abar_j)_{j\in J}$, thereby obtaining a single, countable $R$-module $N^*$, along with a sufficiently generic sequence $\gammabar=(\gamma_{j,m}:j\in J,m\in\omega)$.   
Specifically,  take $M=RZ$, let $H^*= (h_{j,m,a}:j\in J,m\in\omega,a\in A_j)$ be a fixed realization of $q^{\Abar}$, 
 and  let $N^*$ be a fixed f.g.\ hull over $MH^*$.   

The salient feature of $N^*$ is its omitting types condition. 
As $N^*$ is isomorphic to $M(\Abar,\gammabar)$ from Theorem~\ref{onegamma}, we have
\begin{quotation}  
\noindent For every $j\in J$ and $b\in N^*$, $N^*$ realizes $p_{\gamma_{j,m}}(x,b)$ for every $m\in\omega$ if and only if $b\in A_j+L^{N^*}$.
\end{quotation}

\section{A `mod $L^{RZ}$' Borel reduction}   \label{modLRX}
%
%
%

The goal of this section is to prove Proposition~\ref{almost} that gives a Borel map $\Phi$ from $\K1$ to $\{\hbox{elementary substructures of $N^*$}\}$ that preserves isomorphism going forward, but only recovers the reverse modulo $L^{RX}$.
The map $\Phi$ will be the composition of two maps $F,N$ that we  successively describe below.
We begin by defining suitable substructures of the universal, free-like tagged $R$-module $F^*$.

\begin{Definition}  \label{suitable}
 {\em  A substructure $F'=(RZ',A_j')_{j\in J}$ of $F^*=(RZ,A_j)_{j\in J}$ is {\em suitable} if, for every $j\in J$,
\begin{enumerate}
\item  $A'_j=A_j\cap RZ'$ (this is required to be a substructure);
\item  $A'_*=A_*=RX$; and 
\item  There is a basis $B_j$ of the underlying free module $RZ$ for which
\begin{itemize}
\item  $B_j\cap RZ'$ is a basis for $RZ'$;
\item  $B_j\cap A_j$ is a basis for $A_j$; and
\item  $B_j\cap A_j'$ is a basis for $A_j'$.
\end{itemize}
\end{enumerate}
}
\end{Definition}

Note that if $F'$ is a suitable substructure of $F^*$, then its universe $RZ'$ is an infinite dimensional direct summand of $RZ$, hence $RZ'\preceq RZ\preceq\C$ as $R$-modules.
Let 
\begin{quotation}  \noindent $\K1=\{$all tagged $R$-modules $\Vbar=(RX,V_n)_{n\in\omega}$ with universe $RX$ and $L^{RX}\subseteq V_n$ for every $n\in\omega\}$.  
\end{quotation}
Define a function 
$$F:\K1\rightarrow \{\hbox{suitable substructures of $F^*$}\}$$
as follows.   Recall that $F^*$ was defined using bijections $\rho_n:Y_n\rightarrow RX$.

\medskip
\noindent Given $\Vbar=(RX,V_n)_{n\in\omega}\in \K1$,  for each $n\in\omega$ let $Y_n'=\rho_n^{-1}(V_n)$ and 
let $W_n'=\{y'-\rho_n(y'):y'\in Y_n'\}$.  
Let  $Z':=X\cup\bigcup\{Y_n':n\in\omega\}$.  Note that $X\subseteq Z'\subseteq Z$.  Also, for each $n\in\omega$, $RY_n'$ and $RW_n'$ are free $R$-submodules of $RZ'$ with $RY_n'\subseteq RY_n$ and  $RW_n'\subseteq RW_n$.  

The universe of $F(\Vbar)$ is $RZ'$ and, recalling that $J=\{*\}\cup(\omega\times\{0,1\})$, define $A_j'\subseteq A_j$ by:
\begin{itemize}
\item  $A_*'=A_*=RX$;
\item  $A_{(n,0)}'=RY_n'$; and
\item  $A_{(n,1)}'=RW_n'$.
\end{itemize}

It is easily seen that the map $F$ is Borel.    

\begin{Lemma}  \label{Fmap}   For every $\Vbar\in \K1$, $F(\Vbar)$ is a suitable substructure of $F^*$.
\end{Lemma}

\begin{proof}  Showing that $F(\Vbar)$  is a substructure, i.e., verifying that $A_j'=A_j\cap RZ'$ for each $j\in J$ is routine.
Concerning bases, for $j=*$ or $j=(n,0)$, we can simply take $B_j$ to be the set $Z$.  
As $$RZ=RX\oplus \bigoplus\{RY_n:n\in\omega\}$$
we have $RZ'=RX\oplus\bigoplus\{RY_n':n\in\omega\}$, so (2) of Definition~\ref{suitable} is satisfied by $B_j=Z$.
For $j=(n,1)$ we need a new basis for $RZ$ that is suggested by the following Claim.

\medskip
\noindent{\bf Claim.}  For each $n\in \omega$, $RX\oplus RW_n=RX\oplus RY_n$ and 
$RX\oplus RW_n'=RX\oplus RY_n'$.

\begin{proof}  For the first equality, first choose any $w\in W_n$.   Choose $y\in Y_n$ so that $w=y-\rho_n(y)$.
As $y\in RY_n$ and $\rho_n(y)\in RX$, $w\in RX\oplus RY_n$.  It follows immediately that $RW_n\subseteq RX\oplus RY_n$,
hence $RX\oplus RW_n\subseteq RX\oplus RY_n$.  

The other direction is similar.   Choose any $y\in Y_n$.  Let $w=y-\rho_n(y)$, hence $y=w+\rho_n(y)\in RW_n\oplus RX$.
Thus, $RX\oplus RY_n\subseteq RX\oplus RW_n$.  

The argument for $RW_n'$ and $RY_n'$ is identical to the above.
\qed
\end{proof}

Thus, when $j=(n,1)$, note that $$RZ=RX\oplus RY_n\oplus \bigoplus\{RY_k:k\neq n\}$$
so the Claim gives
$$RZ=RX\oplus RW_n\oplus \bigoplus\{RY_k:k\neq n\}$$
Thus, take $B_j=X\cup W_n\cup \bigcup\{Y_k:k\neq n\}$.   Then $B_j$ is visibly a basis for $RZ$, and
verifying Clause (2) of Definition~\ref{suitable} is routine.  
\qed
\end{proof}

Next, we describe a Borel mapping $N:\{$suitable substructures of $F^*\}\rightarrow\{$elementary substructures of $N^*\}$.

Given a suitable substructure $F'=(RZ',A_j')_{j\in J}$ of $F^*$, let $H'=\{h_{j,m,a}:j\in J,m\in\omega,a\in A_j'\}\subseteq H^*$.
Apply Lemma~\ref{Borelpiece} to $N^*$ with $B=(RZ)H^*$ and $A=(RZ')H'$ to get $N(F')\preceq N^*$, $N(F')$ an f.g.\ hull over $(RZ')H'$,
with $\fg {N(F')} {(RZ')H'} {(RZ)H^*}$. As noted there, the map $N$ is Borel.  But note that the pp type of $H'/RZ$ is generated by
pp data whose parameters lie in $RZ'$, so $\fg {H'} {RZ'} {RZ}$.  By transitivity and symmetry of non-forking, the two assertions
imply $\fg {N(F')} {RZ'} {RZ}$.

The following Proposition indicates that the salient feature of $N^*$ passes down to $N(F')$ for any suitable substructure
$F'$ of $F^*$.

\begin{Proposition}  \label{goodness}  Suppose $F'$ is any suitable substructure of $F^*$ and let $N(F')\preceq N^*$ be constructed as above.
Then for every $j\in J$ and $b\in N(F')$,
 $$\hbox{$N(F')$ realizes $p_{\gamma_{j,m}}(x,b)$ for every $m\in\omega \Longleftrightarrow
 b\in A_j'+L^{N(F')}$}$$
 \end{Proposition}

\begin{proof}  First, note that for every $j\in J$ and $m\in\omega$, the element $h_{j,m,a}\in H'$ realizes $p_{\gamma_{j,m}}(x,a)$
for every $a\in A_j'$.  Thus, by Lemma~\ref{typesmodL}, every $b\in A_j'+L^{N(F')}$ has $p_{\gamma_{j,m}}(x,b)$ realized in $N(F')$.

The converse is more interesting.  Fix $j\in J$ and assume $p_{\gamma_{j,m}}(x,b)$ is realized in $N(F')$ for every $m\in\omega$.  Then as $N(F')\preceq N^*$, each of the types is also realized in $N^*$.  Thus, by Theorem~\ref{onegamma},
$b\in A_j+L^{N^*}$.  Choose an element $a\in A_j$ such that $b-a\in L^{N^*}$.   For definiteness, assume $b-a\in L_n^{N^*}$.
By suitability, choose a basis $B_j$ for $RZ$ such that $D:=B_j\cap A_j'$ is a basis for $A_j'$,  $E:=B_j\cap A_j$ is a basis for $A_j$, and $B_j':=B_j\cap RZ'$ is a basis for $RZ'$.  In terms of this basis,
write $$a=\sum r_id_i + \sum s_ke_k$$
where each $d_i\in D$ and each $e_k\in E\setminus D$.  As notation, let $a'=\sum r_id_i$ and let $a^*=\sum s_ke_k$.
We will show that $a^*\in L_n^{N^*}$ by showing that each $e_k$-coefficient $s_k\in L_n^R$.

To see this, fix a coordinate $e_k$ and choose $r\in G_n$ and $m\in\omega$.  
Now, $b-a\in L_n^{N^*}$ means that $\C\models\phi_m(r(b-a))$.  Thus,
$\phi_m(r(x-a))\in\tp^+(b/RZ'a)$.    As $\fg b {RZ'} a$, $\tp(b/RZ'a)$ is generated by $\tp^+(b/RZ')$,
so there is $\theta(x,\dbar)\in \tp^+(b/RZ')$ with $\theta(x,\dbar)\models\phi_m(r(x-a))$.   As $RZ'\preceq\C$, it follows that there is $c\in RZ'$ such that $\C\models \phi_m(r(c-a))$.
As $\phi_m(\C)$ is a subgroup and since $RZ\preceq \C$ with $c,a\in RZ$, 
we conclude that $$RZ\models \phi_m(r(a-c))$$
However, $c\in RZ'$ and $e_k\not\in RZ'$, so the $e_k$-coefficient of $a-c$ is equal to the $e_k$-coefficient of $a$, namely $s_k$.   It follows that we must have $R\models \phi_m(rs_k)$.   As this holds for every $r\in G_n$ and every $m\in\omega$,
we conclude that $s_k\in L_n^R$.    Thus, $a^*\in L_n^{N^*}$.   So $b=a+\ell$ and $a=a'+\ell'$ for some $\ell,\ell'\in L_n^{N^*}.$
As $a'\in A_j'$, this implies $b\in A_j'+L_n^{N^*}$.   However, since $bA_j'\subseteq N(F')\preceq N^*$, it follows by
Fact~\ref{Lpure} that $b\in A_j'+L^{N(F')}$.
\qed
\end{proof}

\begin{Lemma}  \label{def}  Suppose $\Vbar=(RX,V_n)_{n\in\omega}\in \K1$ and $F'=F(\Vbar)\subseteq F^*$, $N(F')\preceq N^*$ are constructed as above.
For each $n\in\omega$, the set $V_n+L^{N(F')}$ is equal to
$$\{a\in RX+L^{N(F')}:\exists b\in RY_n'+ L^{N(F')}\exists c\in RW_n'+L^{N(F')} (c=b-a)\}$$
\end{Lemma}

\begin{proof}   First, choose $x\in V_n+L^{N(F')}$.   Say $x=a+\ell$ with $a\in V_n$ and $\ell\in L^{N(F')}$.
As $Y_n'=\rho_n^{-1}(V_n)$, choose $b\in Y_n'$ with $\rho_n(b)=a$.  Let $c=b-\rho_n(b)\in W_n'$ by definition of $W_n'$.
Let $c^*=c-\ell\in RW_n'+L^{N(F')}$ and $c^*=b-x$.

Conversely, choose $a\in RX+L^{N(F')}$ with $b\in RY_n'+ L^{N(F')}$, $c\in RW_n'+L^{N(F')}$  and $c=b-a$.
Say $a=a_X+\ell_1$, $b=b'+\ell_2$ and $c=c'+\ell_3$ with $a_X\in RX$, $b'\in RY_n'$, $c'\in RW_n'$ and $\ell_1,\ell_2,\ell_3\in L^{N(F')}$.   As $b'\in RY_n'$, let $a'=\rho_n(b')\in V_n$, and let $c^*=b'-a'\in RW_n'$.
As $c=b-a$ we have 
$$c'=b'-a_X+\ell'$$ for some $\ell'\in L^{N(F')}$.   Subtracting $c^*=b'-a'$ yields
$$c'-c^*=a'-a_X+\ell'$$
with $c'-c^*\in RW_n'$ and $a'-a_X\in RX$.
As $RZ'\preceq\C$ and
$$RZ'=RX\oplus RW_n'\oplus\bigoplus_{k\neq n} RY_k'$$
Lemma~\ref{directL} gives $a'-a_X\in L^{RX}$.  Since $a'\in V_n$, $a_X\in V_n+L^{RX}$.  
Thus, $a\in V_n+L^{N(F')}$ as desired.
\qed
\end{proof}

\begin{Proposition} \label{almost}  There is a Borel map $\Phi:\K1\rightarrow \{$elementary substructures of $N^*\}$
such that, for any $\Vbar^1,\Vbar^2\in \K1$ 
\begin{itemize}  
\item If $\Vbar^1\cong \Vbar^2$ as tagged $R$-modules, then $\Phi(\Vbar^1)\cong \Phi(\Vbar^2)$ as $R$-modules; and
\item  If $\Phi(\Vbar^1)\cong \Phi(\Vbar^2)$ as $R$-modules, then there is a tagged $R$-module isomorphism
$$g:(RX/L^{RX},V_n^1/L^{RX})_{n\in\omega}\rightarrow (RX/L^{RX},V_n^2/L^{RX})_{n\in\omega}$$
\end{itemize}
\end{Proposition}

\begin{proof}  Let $\Phi=N\circ F$ be the composition of the two maps described above.  As each of $F$ and $N$ are Borel, so is $\Phi$.  The first half of showing that $\Phi$ is a Borel reduction is easy.  
Suppose  $f:\Vbar^1\cong \Vbar^2$.  For $i=1,2$, write $F(\Vbar^i)=(RZ'_i,A'_{j,i})$ and 
let $H_i'=\{h_{j,m,a}:j\in J,m\in\omega,a\in A'(j,i)\}\subseteq H^*$.

For each $n\in\omega$ and $v\in V_n^1$, let $y_{n,v}:=\rho_n^{-1}(v)$.  
Then $f$ lifts to $f^*:F(\Vbar^1)\cong F(\Vbar^2)$ with $f^*\mr{RX}=f$ via $f^*(y_{n,v})=y_{n,f(v)}$ for each $n\in\omega$ and $v\in V_n^1$.   As $f^*(y_{n,v}-v)=y_{n,f(v)}-f(v)$, we have both $f^*[RY_{n,1}']=RY_{n,2}'$ and  $f^*[RW_{n,1}']=RW_{n,2}'$.
Then $f^*$ induces a bijection $\fhat:H_1'\rightarrow H_2'$ via $\fhat(h_{j,m,a})=h_{j,m,f^*(a)}$.    

The union $f^*\cup \fhat:RZ^1\cup H_1'\rightarrow RZ^2\cup H_2'$ carries the pp partial type generating $\tp(H_1'/RZ^1)$ to the corresponding partial type generating $\tp(H_2'/RZ^2)$, 
hence  $f^*\cup\fhat$ is an elementary bijection by Proposition~\ref{main}.
Thus, $\Phi(\Vbar^1)\cong\Phi(\Vbar^2)$ follows from the uniqueness of f.g.\ hulls given by Proposition~\ref{uniquehull}.

    Conversely, suppose $f:\Phi(\Vbar^1)\rightarrow \Phi(\Vbar^2)$  is an $R$-module isomorphism.  To ease notation, we write $N^1$ and $N^2$ for $\Phi(\Vbar^1)$ and $\Phi(\Vbar^2)$, respectively.  As membership in $L^{N^i}$ ($i=1,2$) is described by sets of pp formulas, it follows that $f[L^{N^1}]=L^{N^2}$, setwise.  As well, for any $b\in N^1$ and any $\gamma\in\F$, $p_\gamma(x,b)$ is realized in $N^1$ if and only if $p_\gamma(x,f(b))$ is realized in $N^2$.   Thus, by Proposition~\ref{goodness}, for every $j\in J$
we have $f[(A_j')^1+L^{N^1}]=(A_j')^2+L^{N^2}$, setwise.   In particular, as $(A_*')^1=(A_*')^2=RX$, we have
$$f[RX+L^{N^1}]=RX+L^{N^2},\ \hbox{setwise}$$
and similarly, $j=(n,0)$ gives $f[(RY_n')^1+L^{N^1}]=(RY_n')^2+L^{N^2}$ and $j=(n,1)$ gives 
$f[(RW_n')^1+L^{N^1}]=(RW_n')^2+L^{N^2}$, both setwise.

Combining this with Lemma~\ref{def} gives, for each $n\in\omega$
$$f[V_n^1+L^{N^1}]=V_n^2+L^{N^2}, \ \hbox{setwise}$$
Noting that each $V_n^i\subseteq RX$ for $i=1,2$, going forward,  we only consider the restriction of $f$ to $RX+L^{N^1}$ (which we still denote as $f$).  
Coupling this with the fact that $f$ maps $L^{N^1}$ onto $L^{N^2}$,  we can quotient out giving an induced $R$-module isomorphism
$$\fbar:(RX+L^{N^1})/L^{N^1} \ \rightarrow \ (RX+L^{N^2})/L^{N^2}$$
that satisfies $\fbar[V_n^1+L^{N^1}/L^{N^1}]=(V_n^2+L^{N^2})/L^{N^2}$, setwise.  
Continuing, by Fact~\ref{Lpure}, $RX\cap L^{N^i}=L^{RX}$ for $i=1,2$.  Furthermore, since our original $\Vbar^1,\Vbar^2\in\K1$,
we have $L^{RX}\subseteq V_n^i$ for all $n\in\omega$ and all $i=1,2$.   From this we obtain $V_n^i\cap L^{N^i}=L^{RX}$ for all $n\in\omega$, $i=1,2$.  
Thus, it follows from the proof of the Second Isomorphism Theorem for modules that there is an  $R$-module isomorphism
$g:RX/L^{RX}\rightarrow RX/L^{RX}$ induced from $\fbar$ that satisfies
$$g[V_n^1/L^{RX}]=V_n^2/L^{RX} \ \hbox{for each $n\in\omega$}$$
Thus, this function $g$ describes a tagged $R$-module isomorphism as requested.
\qed

\end{proof}

\section{Proving Borel completeness}   \label{Borelreduction}

Whereas Proposition~\ref{almost} was a good attempt, in order to get our Borel reduction, we need to handle `modding out by
$L^{RZ}$.'   Crucially, $L^{RZ}=L^R\cdot Z$ by Lemma~\ref{freeL}, so the following idea is somewhat natural.

Recall that by Lemma~\ref{2sided}, $L^R$ is a proper, 2-sided ideal of the ring $R$.  Thus, $S=R/L^R$ is a countable, unital ring.
By Theorem~1.4 of \cite{LU}, the theory of tagged $S$-modules is Borel complete.  Thus, simply by reindexing the tags,
$$\KK=\{\hbox{tagged $S$-modules $(\omega, U_n)_{n\in\omega}$ with $U_0=(0)$}\}$$
is Borel complete.  
Thus, to prove the Borel completeness of $Th(RX)$,
it suffices to find a Borel reduction $\Phi^*:\KK\rightarrow\{$elementary substructures of $N^*\}$.

Note that any $S$-module $M$ can be naturally construed as an $R$-module by defining, for each $r\in R$,
multiplication $r\cdot:M\rightarrow M$ as $r\cdot(a)=(r+L^R)(a)$.  

Fix, once and for all, a bijection $\pi:X\rightarrow\omega$.   

For any $\U=(\omega,U_n)_{n\in\omega}\in\KK$, let 
$\pi_{\U}:RX\rightarrow \U$ be the unique $R$-module homomorphism extending $\pi$.

The following observation illustrates the advantage of working with $R/L^R$-modules.

\medskip
\noindent
{\bf Claim.}  For each $\U\in\KK$, $L^{RX}=L^R\cdot X\subseteq\ker(\pi_{\U})$.

\begin{proof}  The first equality is Lemma~\ref{freeL}.  For the second, choose any $a\in L^{RX}$.  
By the first equality,  $a=\sum r_ix_i$ with each $r_i\in L^R$ and each $x_i\in X$. 
Then, as $\pi_{\U}$ is a homomorphism, our definition of $r\cdot$ gives
$\pi_{\U}(a)=\sum \pi_{\U}(r_ix_i)=0$, so $a\in \ker(\pi_{\U})$.
\end{proof}

Given this, we define a Borel map $h:\KK\rightarrow\K1$ as follows.
Given a tagged $S$-module $\U=(\omega,U_n)_{n\in\omega}$ with $U_0=(0)$,  for each $n\in\omega$,
let $V_n=\pi_{\U}^{-1}(U_n)$ and define
 $$h(\U):=(RX,V_n)_{n\in\omega}$$
Note that for each $n$, $L^{RX}\subseteq\ker(\pi_{\U})\subseteq V_n$, and that $V_0=\ker(\pi_{\U})$.
Thus, $h(\U)\in \K1$ and
it is easily checked that the map $h:\KK\rightarrow\K1$ is Borel.

Put $\Phi^*:=\Phi\circ h:\KK\rightarrow \{$elementary substructures of $N^*\}$.
Then, combining our comments above with Proposition~\ref{almost}, $\Phi^*$ is Borel.
We show that $\Phi^*$ preserves isomorphism in both directions.  Fix $\U_1,\U_2\in \KK$ and let
$\pi_1,\pi_2$ be the maps $\pi_{\U_1},\pi_{\U_2}$ described above.  Note that as $U_0^1=U_0^2=0$,
$V_0^i=\ker(\pi_i)$ for $i=1,2$.  

First, suppose $f:\U_1\rightarrow\U_2$ is a tagged $S$-module isomorphism.  
Then as $f$ is a permutation of $\omega$, it induces a permutation of $X$, and hence an $R$-module
automorphism $f^*:RX\rightarrow RX$. A diagram chase gives that $\pi_{2}f^*=f\pi_{1}$, hence
$$f^*[\pi_{1}^{-1}(U_n^1)]=\pi_{2}^{-1}(U^2_n)$$
for each $n\in\omega$.  Thus, $f^*:h(\U_1)\rightarrow h(\U_2)$ is a tagged $R$-module isomorphism.  
By Proposition~\ref{almost}, we conclude that $\Phi^*(\U_1)\cong \Phi^*(\U_2)$.

Conversely, suppose $\Phi^*(\U_1)\cong \Phi^*(\U_2)$ as $R$-modules.  
By Proposition~\ref{almost} again, there is a tagged $R$-module isomorphism
$$g:(RX/L^{RX},V_n^1/L^{RX})_{n\in\omega}\rightarrow (RX/L^{RX},V_n^2/L^{RX})_{n\in\omega}$$
In particular, the tag-0 condition gives 
$$g[\ker(\pi_1)/L^{RX}]=\ker(\pi_2)/L^{RX}$$
Thus, $g$ induces an isomorphism $g^*:RX/\ker(\pi_1)\rightarrow RX/\ker(\pi_2)$ with
$g^*[V_n^1/\ker(\pi_1)]=V_n^2/\ker(\pi_2)$ for each $n\in\omega$.
Also, via the identifications given by the First Isomorphism Theorem, $RX/\ker(\pi_i)\cong \U_i$ for  $i=1,2$.
Composing these maps gives the required tagged $S$-module isomorphism.
\qed
\end{proof}

\section{Proofs of Corollaries~\ref{BassCor} and \ref{simplecase}}  \label{last}

\noindent{\bf Proof of Corollary~\ref{BassCor}.}

\medskip

By Bass's Theorem P from \cite{Bass},  the ring $R$ is left perfect if and only if
the principal right ideals of $R$ satisfy the descending chain
condition.  Hence, if $R$ is not left perfect, there are
$(a_n:n\in \omega)$ from $R$ such that
$$a_0R\supsetneq a_1R\supsetneq a_2R\supsetneq\cdots.$$
For each $n$, the principal right ideal $a_nR$ is the solution
set in $\RsubR$ of the pp formula $\varphi_n(x):=\exists y\,(a_ny=x)$.
Thus,  if a countable ring $R$ is not left perfect, then $R$ satisfies the assumptions of Theorem~\ref{big2}, so $T_R$ is Borel complete.

\medskip

\noindent{\bf Proof of Corollary~\ref{simplecase}.}

 $(1)\Rightarrow(2)$:   Assume $R$ is a countable, simple ring such that $T_R$ is not Borel complete.  By Corollary~\ref{BassCor}, $R$ is left perfect.  Thus, by \cite{Bass}, $R$ is semilocal, i.e., $R/J(R)$ is semisimple artinian, where
$J(R)$ denotes the Jacobson radical.  However, as $J(R)$ is always a 2-sided ideal, $R$ being simple implies $J(R)=0$, hence $R$ is semisimple artinian.   As $R$ was assumed to be simple, it is simple artinian.

$(2)\Rightarrow(3)$: This is immediate by the Artin-Wedderburn classification of semisimple artinian rings.


$(3)\Rightarrow(1)$:  Let $e=e_{11}$, the $n\times n$  diagonal matrix with all zero entries except a 1 in the $(1,1)$ position.
By Morita equivalence, the functor $M\mapsto eM$ is an equivalence between left $M_n(D)$-modules and left $D$-vector spaces, with quasi-inverse $V\mapsto M_n(D)e\otimes_D V$.   Thus, a left $M_n(D)$-module $M$ is determined up to isomorphism by the cardinal $\dim_D(eM)$.   Since $D$ is countable, there are only countably many possibilities for $\dim_D(V)$ for any countable vector space $V$, hence there are only countably many isomorphism classes of countable left $M_n(D)$-modules.  In particular, $T_R$ is not Borel complete.
\qed

 \end{document}